\documentclass[12pt]{amsart}
\usepackage[final]{graphicx}
\usepackage{subcaption}
\usepackage{float}
\usepackage[colorlinks=true]{hyperref}
\usepackage{accents}
\usepackage{color}
\newcommand{\beq}{\begin{eqnarray*}}
\newcommand{\feq}{\end{eqnarray*}}
\newcommand{\beqn}{\begin{eqnarray}}
\newcommand{\feqn}{\end{eqnarray}}

\newtheorem{theorem}{Theorem}[section]
\newtheorem{lemma}[theorem]{Lemma}

\theoremstyle{definition}

\newtheorem{example}[theorem]{Example}

\theoremstyle{remark}

\numberwithin{equation}{section}

\newtheorem*{theorem*}{Theorem}

\begin{document}
\title[Traffic flows with looking ahead$\&$behind dynamics]{Traffic flow models with looking ahead-behind dynamics}

\author{Yongki Lee$^\dag$}
\address{$^\dag$Department of Mathematical Sciences, Georgia Southern University, Statesboro,  Georgia 30458}

\email{yongkilee@georgiasouthern.edu}
\begin{abstract} 
Motivated by the traffic flow model with Arrhenius look-ahead relaxation dynamics introduced in [A. Sopasakis and M.A. Katsoulakis, SIAM J. Appl. Math., 66 (2006), p. 921--944], this paper proposes a traffic flow model with look ahead relaxation-behind intensification by inserting look behind intensification dynamics to the flux. Finite time shock formation conditions in the proposed model with various types of interaction potentials are identified. Several numerical experiments are performed in order to demonstrate the performance of the modified model. It is observed that, compare to other well-known traffic flow models, the model equipped with look ahead relaxation-behind intensification has both enhanced dispersive and smoothing effects.\\

\end{abstract}
\maketitle

\section{Introduction}
In this paper, we are concerned with the shock formation phenomena - bounded solutions with unbounded derivatives -  for a class of nonlocal conservation laws,
\begin{equation}\label{1main}
\left\{
  \begin{array}{ll}
    \partial_t u + \partial_x F(u, \bar{u}) =0, & t>0, x \in \mathbb{R}, \\
    u(0,x)=u_0 (x), &   x\in \mathbb{R},\hbox{}
  \end{array}
\right.
\end{equation}
where $u$ is the unknown, $F$ is a given smooth function, and $\bar{u}$ is given by
\begin{equation}\label{ubar}
\bar{u}(t,x)=(K*u)(t,x)=\int_{\mathbb{R}} K(x-y)u(t,y) \, dy,
\end{equation}
where $K$ to be chosen later. The nonlinear advection couples both local and nonlocal mechanism. 

This class of conservation laws, identified in \cite{LL15}, appears in several applications including traffic flows \cite{SK06, KP09}, the collective motion of biological cells \cite{MDS08, DS05}, dispersive water waves \cite{GW74, Holm, Liu0}, 
high-frequency waves in relaxing medium \cite{Hunter, Vak1} and the kinematic sedimentation model \cite{Kynch, KZ99, BBKT}.

There are some distinguished special cases of \eqref{1main} with the kernel $K$:\\
$\bullet$ A shallow water model proposed by Whitham \cite{GW74}
$$u_t + \frac{3c_0}{2h_0}uu_x + \bar{u}_x =0,$$
corresponding to \eqref{1main} with $F(u, \bar{u})=\frac{3c_0}{4h_0} u^2 + \bar{u}$ and $K(r)=\frac{\pi}{4}\exp(-\pi |r|/2)$;\\
$\bullet$ The hyperbolic Keller-Segel model with logistic sensitivity \cite{DS05}
\begin{equation}
\left\{
  \begin{array}{ll}
     u_t + [u(1-u)\partial_x S]_x=0,\\
    - S_{xx} + S = u, & \hbox{}
  \end{array}
\right.
\end{equation}
corresponding to \eqref{1main} with $F(u, \bar{u})=u(1-u)\bar{u}$ and $K(r)=\partial_r (e^{-|r|}/2)$;\\
$\bullet$ A nonlocal dispersive equation modeling particle suspensions \cite{JR90, RK89, KZ99}
$$u_t + u_x + ((K_a * u)u)_x =0,$$
corresponding to \eqref{1main} with $F(u, \bar{u})=u + \bar{u}u$, $K_a (r)=a^{-1} K(r/a)$ and
\begin{equation}
K(r)=
\left\{
  \begin{array}{ll}
  2/(3(r^2 /4 -1)) & $if$ \ |r| < 2, \\
0, &   $otherwise$.\hbox{}
  \end{array}
\right.
\end{equation}

Along with the above nonlocal models, the model that motivates the present work is the traffic flow model with looking ahead relaxation, introduced by Sopasakis and Katsoulakis:\\
$\bullet$ A traffic flow model with Arrhenius look-ahead dynamics \cite{SK06}
\begin{equation}\label{trafficA}
u_t + [u(1-u)e^{-\bar{u}}]_x =0,
\end{equation}
corresponding to \eqref{1main} with $F(u, \bar{u})=u(1-u)e^{-\bar{u}}$ and 
\begin{equation}\label{traffic_vanilla_kernel}
K(r)=
\left\{
  \begin{array}{ll}
K_0 / \gamma_a, & $if$ -\gamma_a \leq r \leq 0, \\
0, &   $otherwise$;\hbox{}
  \end{array}
\right.
\end{equation}
Here, $u(t,x)$ represents a vehicle density normalized in the interval $[0,1]$, $\gamma_a$ is a positive constant proportional to the look-ahead distance and $K_0$ is a positive interaction strength. This model takes into account interactions of every vehicle with other vehicles ahead within the look ahead distance $\gamma_a$.

Some careful numerical study of the above traffic flow model is carried out in \cite{KP09}. In addition to this, an improved interaction potential 
\begin{equation}\label{linear_kernel}
K(r):=
\left\{
  \begin{array}{ll}
\frac{2K_0}{\gamma_a}(1+\frac{r}{\gamma_a}), & $if$ -\gamma_a \leq r \leq 0, \\
0, &   $otherwise$,\hbox{}
  \end{array}
\right.
\end{equation}
is introduced in \cite{KP09}. This linear interaction potential is intended to take into account the fact that a vehicle's speed is affected more by nearby vehicles than distant ones. In the case of a good visibility(large $\gamma_a$), the numerical examples in \cite{KP09} suggest that \eqref{trafficA} with the linear potential yields solutions that seem to better correspond to reality.

Several finite time shock formation scenarios of solution to  \eqref{trafficA} with \eqref{traffic_vanilla_kernel} were presented in \cite{LL11}. The authors in \cite{LL15}  identified threshold conditions for the finite time shock formation of \eqref{trafficA} subject to two different potentials above. The sub-thresholds for finite time shock formation conditions in \cite{LL15} are consistent with the numerical result in \cite{KP09}. 

We set $K_0=1$ through out this paper, since this parameter is not essential in our blow-up analysis. Then we can rewrite nonlocal term $\bar{u} = K*u$ associated with \eqref{traffic_vanilla_kernel} and \eqref{linear_kernel} as follows, respectively:
\begin{equation}\label{a_const}
\bar{u}(t,x)=\frac{1}{\gamma_a}\int^{x+\gamma_a} _{x} u(t,y) \, dy,
\end{equation}
and
\begin{equation}\label{a_lin}
\bar{u}(t,x)=\frac{2}{\gamma_a}\int^{x+\gamma_a} _{x}  \bigg{(} 1+ \frac{x-y}{\gamma_a}  \bigg{)} u(t,y) \, dy.
\end{equation}


In this paper, we extend \eqref{trafficA} by considering a look behind intensification. That is, we consider the traffic flow model with look ahead relaxation and \emph{look behind intensification}:
\begin{equation}\label{trafficAB}
\left\{
  \begin{array}{ll}
    \partial_t u + \partial_x (u(1-u)e^{-\bar{u} + \tilde{u} }) =0, & t>0, x \in \mathbb{R}, \\
    u(0,x)=u_0 (x), &   x\in \mathbb{R}.\hbox{}
  \end{array}
\right.
\end{equation}
Here, $\bar{u}$ is given in \eqref{a_const} and \eqref{a_lin}. For the $\tilde{u}=K_b * u$, we shall consider constant and linear interaction potentials. i.e.,
\begin{equation}\label{look_b_const}
K_b(r)=
\left\{
  \begin{array}{ll}
K_0 / \gamma_b, & $if$ \ 0 \leq r \leq \gamma_b, \\
0, &   $otherwise$;\hbox{}
  \end{array}
\right.
\end{equation}
and
\begin{equation}\label{look_b_linear}
K_b(r):=
\left\{
  \begin{array}{ll}
\frac{2K_0}{\gamma_a}(1 -\frac{r}{\gamma_b}), & $if$  \ 0  \leq r \leq  \gamma_b, \\
0, &   $otherwise$.\hbox{}
  \end{array}
\right.
\end{equation}
Here, $\gamma_b$ is a nonnegative constant proportional  to the  the look-behind distance. By setting $K_0=1$ again, we can rewrite the nonlocal term $\tilde{u}$ associated with the above kernels as follows, respectively:
\begin{equation}\label{b_const}
\tilde{u}(t,x):=\frac{1}{\gamma_b}\int^{x} _{x-\gamma_b} u(t,y) \, dy,
\end{equation}
and
\begin{equation}\label{b_lin}
\tilde{u}(t,x)=\frac{2}{\gamma_b}\int^{x} _{x-\gamma_b}  \bigg{(} 1- \frac{x-y}{\gamma_b}  \bigg{)} u(t,y) \, dy.
\end{equation}

This look behind intensification model is intended to take into account the driving behavior of some drivers who actively react to vehicle distributions of one's ahead and behind. More precisely, considering the flux in \eqref{trafficAB}, a vehicle's velocity is determined by $(1-u) e^{-\bar{u} +\tilde{u}}$.  Here, the \emph{local} traffic density $u$ at one's location plays an major roles.  In addition to this, \emph{averaged} ahead and behind traffic densities of each driver adjust the velocity via relaxation and intensification effects. The drivers equipped with this strategy monitor/compare the densities ahead and behind within assigned distances, and prefer to choose accelerate(decelerate) when one has a relatively high density behind(ahead). It is natural to assume that $\gamma_a \geq \gamma_b >0$.

The motivation of the behind intensification may be explained as follows.  On one lane highway, if there is a relatively high density behind a driver, the driver may feel that he/she is holding up the traffic and want to i) escape from pressure by driving faster if conditions allow, ii) accelerate to improve the flow of entire convoy of cars.



The objective of this article is twofold:

i) We identify threshold conditions for the finite time shock formation of the traffic flow model \eqref{trafficAB} subject to two different potentials \eqref{a_const}-\eqref{b_const} and \eqref{a_lin}-\eqref{b_lin}. As is known that the typical well-posedness result(e.g., \cite{LL15}) asserts that either a solution of a hyperbolic partial differential equation exists for all time or else there is a finite time such that some norm of the solution becomes unbounded as the life span is approached. The natural question is whether there is a critical threshold for the initial data such that the persistence of the solution regularity depends only on crossing such a critical threshold. This concept of critical threshold and associated methodology is originated and developed in a series of papers by Engelberg, Liu and Tadmor \cite{ELT01, LT02,  TT14} for a class of Euler-Poisson equations.

ii) We investigate performance of the proposed model \eqref{trafficAB} via numerical examples in comparison with \eqref{trafficA} and the Lighthill - Whitham - Richards(LWR) \cite{LW55, R56} model 
\begin{equation}\label{LWR0}
\partial_t u + \partial_x (u(1-u)) =0.
\end{equation}
We are interested in both \emph{dispersive} and \emph{smoothing} effects of the proposed global flux. It is well known that the $u_x$ of the LWR model blows-up if there is a point such that $u_x(t, x) \big{|}_{t=0} >0$. Indeed, the derivative of $d=u_x$ along the characteristic satisfies the Riccatti equation
$$\dot{d} = d^2$$
that leads to the blow-up of $d$ \emph{unless} $d\big{|}_{t=0} \leq 0$.  For some circumstances, this result is unrealistic because no shock formation is observed in the free flow(i.e., sparse traffic). 

On the other hand, in \cite{KP09}, it is observed the global flux in look-ahead model \eqref{trafficA} has some smoothing effect that seems to be able to prevent the shock formation. But simultaneously, the relaxation makes waves less dispersive. In our numerical examples in Section \ref{section3}, it is observed that the proposed model \eqref{trafficAB} greatly \emph{improves} waves' dispersion and smoothing phenomena.


Now, the finite time shock formation results are collectively stated as follows. In the theorems, we assume that $\gamma_a \geq \gamma_b >0$.

\begin{theorem}\label{main1}
(Constant interaction potentials) Consider \eqref{trafficAB} with \eqref{a_const} and \eqref{b_const}. Suppose that $u_0 \in H^2$ and $0\leq u_0 (x) \leq 1$ for all $x\in \mathbb{R}$. If
\begin{equation}\label{const_AB_blowup}
\sup_{x\in \mathbb{R}}[u' _0 (x)] > \frac{\gamma_a + \gamma_b}{\gamma_a \gamma_b} \bigg{(}  \frac{1}{2} + \frac{\sqrt{2}}{4}\sqrt{3 - \min \bigg{\{}-1, \frac{\gamma_a \gamma_b}{\gamma_a + \gamma_b} \inf_{x\in \mathbb{R}} [u'_0 (x)]  \bigg{ \}} } \bigg{)}, 
\end{equation}
then $u_x$ must blow-up at some finite time.
\end{theorem}

\begin{theorem}\label{main2}
(Linear interaction potentials) Consider \eqref{trafficAB} with \eqref{a_lin} and \eqref{b_lin}. Suppose that $u_0 \in H^2$ and $0\leq u_0 (x) \leq 1$ for all $x\in \mathbb{R}$. If there is a point $x \in \mathbb{R}$ such that
\begin{equation}\label{lin_AB_blowup}
u' _0 (x) > \frac{\gamma_a + \gamma_b}{\gamma_a \gamma_b} \bigg{(}  1 + \sqrt{\frac{3}{2} + \big{ ( } \frac{\gamma_a}{2(\gamma_a + \gamma_b)} \big{)}^2} \bigg{)}, 
\end{equation}
then $u_x$ must blow-up at some finite time.
\end{theorem}

Regarding these results some remarks are in order.

i) The condition \eqref{const_AB_blowup} reflects some balance between $\sup_{x\in \mathbb{R}}[u' _0 (x)]$ and $\inf_{x\in \mathbb{R}}[u' _0 (x)]$. It seems the nonpositive term $\inf_{x\in \mathbb{R}}[u' _0 (x)]$ is more negative, then $\sup_{x\in \mathbb{R}}[u' _0 (x)]$ needs to be large for the finite time shock formation.  It can be interpreted that not only the car density behind the traffic jam but also the car density ahead of the traffic jam contribute to the formation of shock.

ii) There is no direct comparison between  \eqref{const_AB_blowup} and look-ahead only model's(\eqref{trafficA}-\eqref{a_const}) blow-up condition
\begin{equation}\label{const_A_blowup}
\sup_{x\in \mathbb{R}}[u' _0 (x)] > \frac{1}{\gamma_a } \bigg{(}  \frac{1}{2} + \frac{\sqrt{2}}{4}\sqrt{3 - \min \bigg{\{}-1, \gamma_a \inf_{x\in \mathbb{R}} [u'_0 (x)]  \bigg{ \}} } \bigg{)},
\end{equation}
which is obtained in \cite{LL15}. Roughly speaking, if $\inf_{x\in \mathbb{R}} [u'_0 (x)]$ is not too negative, the threshold in \eqref{const_AB_blowup} is higher than the one in \eqref{const_A_blowup}. Also, for example, if $\gamma_a =1$ and $\gamma_b =0.5$, then the right hand side of  \eqref{const_AB_blowup} is bigger than that of  \eqref{const_A_blowup} for \emph{any} initial data $u_0 (x)$. The same thing holds for majorities of $\gamma_a$ and $\gamma_b$ values. This is interesting because, even though we insert behind intensification $\tilde{u}$  which in turn increases the waves speed, the blow-up thresholds actually increased. These blow-up results are consistent with numerical experiments in Section \ref{section3}.

iii) In contrast to \eqref{const_AB_blowup}, the blow-up condition in \eqref{lin_AB_blowup} depends only on $\gamma_a$, $\gamma_b$ and the size of the initial slope.  It is interesting to observe that a little difference in the interaction kernels leads to different types of blow-up conditions; one involves \emph{two global terms}  $\sup_{x\in \mathbb{R}}[u' _0 (x)]$ and $\inf_{x\in \mathbb{R}}[u' _0 (x)]$, the other involves \emph{local} $u'(x)$ term only. Furthermore, the blow-up threshold in \eqref{lin_AB_blowup} is higher than the one in \eqref{const_AB_blowup} for not too negative $\inf_{x\in \mathbb{R}} [u'_0 (x)]$.
Indeed, Example \ref{eg4} in Section \ref{section3} shows that the model with linear interaction potentials has less steep wave(behind the traffic jam) than the model with constant potentials for all time.

We now summarize the main arguments in our proofs. We want identify some threshold condition for the shock formation of solutions to \eqref{trafficAB}. The local existence and blow-up alternative results for \eqref{1main} or \eqref{trafficAB} can be found in \cite{LL15}. From Corollary 1 in \cite{LL15}, it suffices to track the dynamics of $u_x$.  The  idea is based on tracing $M(t):=\sup_{x\in \mathbb{R}}[u_x (x,t)]$ and $N(t):=\inf_{x\in \mathbb{R}}[u_x (x,t)]$. The existence and differentiability(in almost everywhere sense) of $M(t)$ and $N(t)$ are proved by Constantin and Escher \cite{CE98}, which we summarize in the following.

\begin{lemma}\label{lem1}(Theorem 2.1 in \cite{CE98})
Let $T>0$ and $u \in C^1 ([0,T];H^2)$. Then for every $t \in [0,T]$ there exists at least one point $\eta(t) \in \mathbb{R}$ with
$$N(t):=\inf_{x \in \mathbb{R}} [u_x (t,x)]=u_x (t, \eta(t)),$$
and the function $N$ is almost everywhere differentiable on $(0,T)$ with
$$
\frac{dN}{dt}(t)=u_{tx}(t, \eta(t)) \ \ a.e. \ \ on \ (0,T).
$$
\end{lemma}

We also state a useful result, which is obtained in \cite{LL09}.
\begin{lemma} (Lemma 3.3 in \cite{LL09})\label{lem2}
Consider the following differential inequality,
\begin{equation}\label{liu_inequality2}
\frac{dB}{dt} \geq a(t)(B - b_1 (t))(B - b_2(t)), \ \ \ B(0)=B_0,
\end{equation}
with $a(t)>0$, $b_1 (t) \leq b_2 (t)$ and that $a(t)$, $b_1 (t)$, $b_2 (t)$ are uniformly bounded.\\
i) If $B_0 > \max b_2$, then $B(t)$ will experience a finite time blow-up.\\
ii) $\min\{B_0 , \min b_1 \} \leq B(t)$, for $t \geq 0$ as long as $B(t)$ remains finite on the time interval $[0,t]$.
\end{lemma}
We remark that the above lemma remains valid even if \eqref{liu_inequality2} holds almost everywhere.

We now conclude this section by outlining the rest of the paper. In Section \ref{section2}, we prove Theorems \ref{main1} and \ref{main2}. In Section \ref{section3}, we demonstrate the performance of the proposed model via several numerical experiments.

\section{Proof of theorems}\label{section2}

\textbf{Proof of Theorem \ref{main1}:}
Let $d:=u_x$ and apply $\partial_x$ to the first equation of \eqref{trafficAB},
\begin{equation}\label{d_eqn}
\begin{split}
\dot{d}&:=(\partial_t +(1-2u)e^{-\bar{u} + \tilde{u}}\partial_x)d\\
&=e^{-\bar{u}+\tilde{u}}\bigg{[} 2d^2 +2(1-2u)(\bar{u}_x -\tilde{u}_x )d  -u(1-u)\big{\{}(-\bar{u}_x  + \tilde{u}_x)^2 + (-\bar{u}_{xx} + \tilde{u}_{xx}   )\big{\}} \bigg{]}.
\end{split}
\end{equation}
Define for $t \in [0,T)$,
\begin{equation}
\begin{split}
&M(t):=\sup_{x \in \mathbb{R}}[u_x (t,x)]=d(t, \xi(t)),\\
&N(t):=\inf_{x \in \mathbb{R}}[u_x (t,x)]=d(t, \eta(t)).\\
\end{split}
\end{equation}

Then, along $(t,\xi(t))$, using \eqref{a_const} and \eqref{b_const}, we have
\begin{equation}\label{const_d_2nd}
\begin{split}
-\bar{u}_{xx} + \tilde{u}_{xx} &= -\frac{1}{\gamma_a} u_x (\xi + \gamma_a) + \bigg{(}  \frac{1}{\gamma_a} + \frac{1}{\gamma_b}   \bigg{)}u_x (\xi) -\frac{1}{\gamma_b}u_x (\xi - \gamma_b),  \\
&\leq  \bigg{(}  \frac{1}{\gamma_a} + \frac{1}{\gamma_b}   \bigg{)}(M -N),
\end{split}
\end{equation}
and
\eqref{d_eqn} can be written as,
\begin{equation}\label{eqn_M}
\begin{split}
\dot{M} &=  e^{-\bar{u}+\tilde{u}}\bigg{[} 2M^2 +2(1-2u)(\bar{u}_x -\tilde{u}_x )M  -u(1-u)\big{\{}(-\bar{u}_x  + \tilde{u}_x)^2 + (-\bar{u}_{xx} + \tilde{u}_{xx}   )\big{\}} \bigg{]} \quad a.e.\\
& \geq  e^{-\bar{u}+\tilde{u}}\bigg{[} 2M^2 +2(1-2u)(\bar{u}_x -\tilde{u}_x )M \\ 
& \ \ \ \ \ \ \ \ \ \  \ \ \ \ \ \ \ \ \ \ \ \ \ \ \ \ \ \ \ \ \ \ \ \ \ \ \ \ -u(1-u)\bigg{\{}(-\bar{u}_x  + \tilde{u}_x)^2 + \bigg{(}  \frac{1}{\gamma_a} + \frac{1}{\gamma_b}   \bigg{)}(M -N)\bigg{\}} \bigg{]} \quad a.e.
\end{split}
\end{equation}

Along $(t, \eta(t))$,  we have
$$-\bar{u}_{xx} +\tilde{u}_{xx} =  -\frac{1}{\gamma_a} u_x (\eta  + \gamma_a) + \bigg{(}  \frac{1}{\gamma_a} +\frac{1}{\gamma_b} \bigg{)}N - \frac{1}{\gamma_b}u_x (\eta - \gamma_b)  \leq 0 ,$$
and
\eqref{d_eqn} can be written as,
\begin{equation}\label{eqn_N}
\begin{split}
\dot{N} &= e^{-\bar{u}+\tilde{u}}\bigg{[} 2N^2 +2(1-2u)(\bar{u}_x -\tilde{u}_x )N  -u(1-u)\big{\{}(-\bar{u}_x  + \tilde{u}_x)^2 + (-\bar{u}_{xx} + \tilde{u}_{xx}   )\big{\}} \bigg{]} \quad a.e. \\
&\geq e^{-\bar{u}+\tilde{u}}\bigg{[} 2N^2 +2(1-2u)(\bar{u}_x -\tilde{u}_x )N  -u(1-u)(-\bar{u}_x  + \tilde{u}_x)^2 \bigg{]} \quad a.e.
\end{split}
\end{equation}
Now, we write the inequality in \eqref{eqn_N} as
\begin{equation}\label{Nt}
\dot{N} \geq 2 e^{-\bar{u} + \tilde{u} }(N - N_1)(N -N_2) \quad a.e.\;,
\end{equation}
where $$
N_1 (u , \bar{u}_x, \tilde{u}_x) = \frac{-(1-2u)(\bar{u}_x -\tilde{u}_x)  - \sqrt{\{(1-2u)(\bar{u}_x - \tilde{u}_x ) \}^2 + 2u(1-u)(-\bar{u}_x + \tilde{u}_x )^2  } }{2}
$$ and
$$
N_2 (u, \bar{u}_x, \tilde{u}_x) =\frac{-(1-2u)(\bar{u}_x -\tilde{u}_x)  + \sqrt{\{(1-2u)(\bar{u}_x - \tilde{u}_x ) \}^2 + 2u(1-u)(-\bar{u}_x + \tilde{u}_x )^2  } }{2}.
$$
We note that $N_1 \leq 0 \leq N_2$ because $0 \leq u(t) \leq 1$.
It can be shown later that $N_1$  is uniformly bounded from below,
\begin{equation}\label{N1l}
N_1 \geq -\frac{\gamma_a + \gamma_b}{\gamma_a \gamma_b}.
\end{equation}
Applying Lemma \ref{lem2} part (ii) to (\ref{Nt}) with $\min_{0 \leq u \leq 1 , \ |\omega_i|\leq \frac{1}{\gamma}} N_1 (u, \omega_1 , \omega_2)=-\frac{\gamma_a + \gamma_b}{\gamma_a \gamma_b}$, we obtain
\begin{equation}\nonumber
N(t) \geq \min\bigg{\{} -\frac{\gamma_a + \gamma_b}{\gamma_a \gamma_b}, \ N(0)  \bigg{\}}=:   \bigg{(}  \frac{\gamma_a + \gamma_b}{\gamma_a \gamma_b} \bigg{)} \tilde{N}_0 .
\end{equation}

Substituting this lower bound into \eqref{eqn_M}, we obtain
\begin{equation*}
\begin{split}
\dot{M} &\geq e^{-\bar{u}+\tilde{u}} \bigg{(} 2M^2 + \bigg{\{} 2(1-2u)(\bar{u}_x - \tilde{u}_x) -\frac{u(1-u)(\gamma_a + \gamma_b)}{\gamma_a \gamma_b} \bigg{\}} M\\
& -u(1-u) (-\bar{u}_x + \tilde{u}_x)^2  + \frac{u(1-u)(\gamma_a +\gamma_b)^2\tilde{N}_0}{(\gamma_a \gamma_b)^2} \bigg{)}
\ a.e.\;
\end{split}
\end{equation*}
Rewriting of this inequality gives
\begin{equation}\label{eqn_M_1}
\dot{M} \geq 2 e^{-\bar{u} + \tilde{u}}(M-M_1)(M -M_2)\; \ \  a.e.\;,
\end{equation}
where $M_2(\geq M_1)$ is given by
\begin{equation*}
\begin{split}
M_2 &:= \frac{-\{2(1-2u)(\bar{u}_x -\tilde{u}_x)-\frac{u(1-u)(\gamma_a + \gamma_b)}{\gamma_a \gamma_b} \}}{4}\\
&+\frac{\sqrt{\{2(1-2u)(\bar{u}_x -\tilde{u}_x)-\frac{u(1-u)(\gamma_a + \gamma_b)}{\gamma_a \gamma_b} \} ^2  +8u(1-u)(-\bar{u}_x + \tilde{u}_x)^2 -\frac{8u(1-u)(\gamma_a + \gamma_b)^2 \tilde{N}_0}{(\gamma_a \gamma_b)^2}}}{4}
\end{split}
\end{equation*}

We claim that $M_2$ has an uniform upper bound,
\begin{equation}\label{M2}
    M_2 \leq \frac{\gamma_a + \gamma_b}{\gamma_a \gamma_b} \bigg{[} \frac{1}{2}+ \frac{\sqrt{2}}{4}\cdot \sqrt{3-\tilde{N}_0}  \bigg{]}.
\end{equation}
By Lemma \ref{lem2} (i), if
$$
M(0)>\frac{\gamma_a + \gamma_b}{\gamma_a \gamma_b} \bigg{[} \frac{1}{2}+ \frac{\sqrt{2}}{4}\cdot \sqrt{3-\tilde{N}_0}  \bigg{]},
$$
then $M(t)$ will blow up a finite time. This is exactly the threshold condition as stated in Theorem \ref{main1}.

To complete our proof we still need to verify both claims (\ref{M2}) and (\ref{N1l}).

To verify (\ref{M2}), we set
\begin{equation*}
\begin{split}
v&:=\frac{\gamma_a \gamma_b}{\gamma_a + \gamma_b} (\bar{u}_x - \tilde{u}_x) \\
&=\frac{\gamma_a \gamma_b}{\gamma_a + \gamma_b} \bigg{\{} \frac{1}{\gamma_a}u(x+\gamma_a) - \big{(}\frac{\gamma_a + \gamma_b}{\gamma_a  \gamma_b}\big{)}u(x) + \frac{1}{\gamma_b}u(x-\gamma_b)   \bigg{\}}.
\end{split}
\end{equation*}
From $0\leq u(t) \leq 1$ it follows that $-1 \leq v \leq 1$.  If suffices to find upper bound for $M_2$ over the set
$$
\Omega:=\{(u,v)\in \mathbb{R}^2 \ | \ 0\leq u \leq 1 , \ \ -1 \leq v \leq 1 \}.
$$
In fact,
{\small
\begin{equation}\nonumber
\begin{split}
M_2 &=\frac{-\{2(1-2u)v-u(1-u)\}+\sqrt{\{2(1-2u)v-u(1-u)\}^2+8u(1-u)(v^2 -\tilde{N}_0) } }{4(\gamma_a \gamma_b /(\gamma_a + \gamma_b))}  \\
&\leq \frac{\gamma_a + \gamma_b}{4\gamma_a \gamma_b} \big{[} 2+ \sqrt{ 4 + 2(1 - \tilde{N}_0) } \big{]} . \\
\end{split}
\end{equation}
}Here, we use $\max_{(u,v)\in \Omega} \{ -2(1-2u)v +u(1-u) \}=2$ which can be verified easily since the underlying function is linear in $v$ and quadratic in $u$. For the next one, $\max_{(u,v)\in \Omega} \{ 8u(1-u)(v^2 - \tilde{N}_0) \} = 2(1-\tilde{N}_0)$ is used, which is obtained from the upper bound $u(1-u)\leq 1/4$. Applying Lemma \ref{lem2}, part (i), we obtain the desired result.

Finally, we are left with the verification of \eqref{N1l}. With $v$ defined above, we have
\begin{equation}\nonumber
Q:=\frac{\gamma_a \gamma_b}{\gamma_a + \gamma_b} N_1=\frac{-(1-2u)v - \sqrt{\{(1-2u)v \}^2 +2u(1-u)v^2  } }{2}.
\end{equation}
By rearranging,
\begin{equation}
\begin{split}
Q^2 &= \frac{u(1-u)v^2}{2} -Q\cdot(1-2u)v\\
&\leq \frac{u(1-u)v^2}{2} + \epsilon Q^2 + \frac{(1-2u)^2}{4\epsilon}v^2, \ \ \ 0 < \epsilon < 1.
\end{split}
\end{equation}
It follows that
\begin{equation}
\begin{split}
(1-\epsilon)Q^2  &\leq \frac{v^2}{4 \epsilon}\{ (1-2u)^2 + 2 \epsilon u(1-u) \}\\
&\leq \frac{1}{4 \epsilon},
\end{split}
\end{equation}
where the maximum value is achieved at $\partial \Omega$. This gives
$$
Q^2 \leq \frac{1}{4 \epsilon (1-\epsilon)}.
$$
Since $\epsilon$ is arbitrary, we choose $\epsilon =\frac{1}{2}$ to get $Q^2\leq 1$, hence $Q\geq -1$, which gives (\ref{N1l}). This completes the proof of Theorem \ref{main1}. $\hfill\square$\\

\textbf{Proof of Theorem \ref{main2}:}
Let $d:=u_x$ and apply $\partial_x$ to the first equation of \eqref{trafficAB},
\begin{equation}\label{d_eqn_2}
\begin{split}
\dot{d}&:=(\partial_t +(1-2u)e^{-\bar{u} + \tilde{u}}\partial_x)d\\
&=e^{-\bar{u}+\tilde{u}}\bigg{[} 2d^2 + 2(1-2u)(\bar{u}_x -\tilde{u}_x )d  -u(1-u)\big{\{}(-\bar{u}_x  + \tilde{u}_x)^2 + (-\bar{u}_{xx} + \tilde{u}_{xx}   )\big{\}} \bigg{]}.
\end{split}
\end{equation}

Using \eqref{a_lin} and \eqref{b_lin}, we display the derivatives of the global terms(we omit the $t-$dependence).
\begin{equation}
\bar{u}_x = \frac{2}{\gamma_a} \bigg{[} \frac{1}{\gamma_a} \int^{x+\gamma_a} _{x} u(y) \, dy -u(x)  \bigg{]}, \ \ \bar{u}_{xx} = \frac{2}{\gamma_a} \bigg{[} \frac{1}{\gamma_a} \big{\{} u(x+\gamma_a)  -u(x) \big{\}}  -d(x)   \bigg{]}.
\end{equation}
Also,
\begin{equation}
\tilde{u}_x = \frac{2}{\gamma_b} \bigg{[} -\frac{1}{\gamma_b} \int^{x} _{x-\gamma_b} u(y) \, dy +u(x)  \bigg{]}, \ \ \tilde{u}_{xx} =\frac{2}{\gamma_b} \bigg{[} -\frac{1}{\gamma_b} \big{\{} u(x)  -u(x-\gamma_b) \big{\}}  +d(x)   \bigg{]}.
\end{equation}

We reorder $d(x)$ terms in $ (-\bar{u}_{xx} + \tilde{u}_{xx})$ and obtain
\begin{equation}\label{d_eqn_22}
\begin{split}
\dot{d}=e^{-\bar{u} +\tilde{u}} \bigg{[}&2d^2  +\big{\{}2(1-2u)(\bar{u}_x -\tilde{u}_x ) -2u(1-u)(\frac{1}{\gamma_a} + \frac{1}{\gamma_b})   \big{\}}d\\
 &-u(1-u) \big{\{}  (-\bar{u}_x + \tilde{u}_x)^2  + (\frac{2}{\gamma^2 _a}  -\frac{2}{\gamma^2 _b})u(x)  +\frac{2}{\gamma^2 _b} u(x-\gamma_b) -\frac{2}{\gamma^2 _a} u(x+\gamma_a) \big{\}}   \bigg{]}.
 \end{split}
\end{equation}
We should point out that in contrast to the model with constant potentials(see \eqref{const_d_2nd}), the above equation does \emph{not} have $u_x (x+\gamma_a)$ or $u_x (x -\gamma_b)$(see  \eqref{const_d_2nd}  ). Later this difference leads to the structurally different blow up conditions.

Rewriting of \eqref{d_eqn_22} gives
\begin{equation}
\dot{d}=2e^{-\bar{u} + \tilde{u}} (d-D_1)(d-D_2),
\end{equation}
where $D_2 (\geq D_1)$ is given by
\begin{equation}
D_2 = \frac{-B+\sqrt{B^2 + C}}{2},
\end{equation}
where
$$B=(1-2u)(\bar{u}_x -\tilde{u}_x ) -u(1-u)(\frac{1}{\gamma_a} + \frac{1}{\gamma_b}),$$
and
$$C=2u(1-u) \bigg{\{}  (-\bar{u}_x + \tilde{u}_x)^2  + (\frac{2}{\gamma^2 _a}  -\frac{2}{\gamma^2 _b})u(x)  +\frac{2}{\gamma^2 _b} u(x-\gamma_b) -\frac{2}{\gamma^2 _a} u(x+\gamma_a) \bigg{\}}. 
$$

To find an upper bound of $D_2$, we set 
$$v:= \frac{\gamma_a \gamma_b}{2(\gamma_a + \gamma_b)} (\bar{u}_x - \tilde{u}_x).$$ 
Then, one can see that
\begin{equation}\label{v_bound2}
 | v |  \leq 1.
 \end{equation}
Indeed, since
\begin{equation*}
\begin{split}
v&= \frac{\gamma_a \gamma_b}{2(\gamma_a + \gamma_b)} \bigg{\{ }\frac{2}{\gamma^2 _a} \int^{x+\gamma_a} _{x} u(y) \, dy - \bigg{(} \frac{2}{\gamma_a} +\frac{2}{\gamma_b}  \bigg{)}u(x) +\frac{2}{\gamma^2 _b} \int^x _{x-\gamma_b} u(y) \, dy \bigg{\}},
\end{split}
\end{equation*}
 and $0 \leq u(x) \leq 1$ we obtain the desired bound.

Now consider
\begin{equation}
\begin{split}
 \frac{\gamma_a \gamma_b}{2(\gamma_a + \gamma_b)} D_2 = &\frac{1}{2} \bigg{[} - \big{\{} (1-2u)v - \frac{u(1-u)}{2}   \big{\}}\\
 & + \sqrt{ \big{\{} (1-2u)v - \frac{u(1-u)}{2}   \big{\}}^2  +  \big{(}\frac{\gamma_a \gamma_b}{2(\gamma_a + \gamma_b)} \big{)}^2 C   } \bigg{]},
\end{split}
\end{equation}
where $v$ and $C$ are defined above.

As in the proof of Theorem \ref{main1} we let
$$
\Omega:=\{(u,v)\in \mathbb{R}^2 \ | \ 0\leq u \leq 1 , \ \ -1 \leq v \leq 1 \}.
$$
Then it holds
\begin{equation}\label{gamma_D2}
\begin{split}
 \frac{\gamma_a \gamma_b}{2(\gamma_a + \gamma_b)} D_2 \leq \frac{1}{2} \bigg{[}  1+ \sqrt{1 +   \big{(}\frac{\gamma_a \gamma_b}{2(\gamma_a + \gamma_b)} \big{)}^2 C}  \bigg{]}.
\end{split}
\end{equation}
Here, we use $\max_{(u,v)\in \Omega} \{ (1-2u)v -u(1-u)/2 \}=1$ which can be verified easily since the underlying function is linear in $v$ and quadratic in $u$. 

To find an upper bound of 
  $ \big{(}\frac{\gamma_a \gamma_b}{2(\gamma_a + \gamma_b)} \big{)}^2 C$, we expand and obtain
  \begin{equation*}
\begin{split}
\big{(}\frac{\gamma_a \gamma_b}{2(\gamma_a + \gamma_b)} \big{)}^2 C &= 2u(1-u)\bigg{[} v^2\\   
& \ \ \ +\big{(}\frac{\gamma_a \gamma_b}{2(\gamma_a + \gamma_b)} \big{)}^2  \bigg{\{}  (\frac{2}{\gamma^2 _a}  -\frac{2}{\gamma^2 _b})u(x)  +\frac{2}{\gamma^2 _b} u(x-\gamma_b) -\frac{2}{\gamma^2 _a} u(x+\gamma_a)  \bigg{\}} \bigg{]}\\
&\leq \frac{1}{2} \bigg{[} 1+  \big{(}\frac{\gamma_a \gamma_b}{2(\gamma_a + \gamma_b)} \big{)}^2 \bigg{\{} \frac{2}{\gamma^2 _b}  \bigg{\}}  \bigg{]}.
\end{split}
\end{equation*}
The last inequality holds because of \eqref{v_bound2}, $0\leq u(x) \leq 1$, $u(1-u) \leq 1/4$ and
$$  (\frac{2}{\gamma^2 _a}  -\frac{2}{\gamma^2 _b}) < 0.$$
Therefore,  \eqref{gamma_D2} leads to
$$D_2 \leq  \frac{\gamma_a + \gamma_b}{\gamma_a \gamma_b} \bigg{[} 1+ \sqrt{ \frac{3}{2} + \big{(} \frac{\gamma_a}{2(\gamma_a + \gamma_b)}  \big{)}^2 }\bigg{]}.$$
Applying part (i) of Lemma \ref{lem2}, we complete the proof of Theorem \ref{main2}. $\hfill \square$.

$$$$

\section{Numerical examples.}\label{section3}
In this section, we demonstrate the performance of the proposed look ahead relaxation-behind intensification model \eqref{trafficAB}, in comparison with look ahead relaxation only model \eqref{trafficA} and LWR(non-global flux) model \eqref{LWR0}. 

 In all numerical examples, Lax-Friedrichs numerical scheme is applied, and we take large enough computational domain  so that no waves touch its boundary within the final computational time. The Courant-Friedrichs-Lewy(CFL) numbers are chosen as $0.5$ or $0.25$, and all solutions are computed on uniform grids with $\Delta x =1/100$ or $\Delta x=1/50$. In order to compute the nonlocal term $\bar{u}(t,x)=K*u$, we used the trapezoidal method.

We re-display the models of considerations here for readers' convenience:\\
$\bullet$ (LWR)  Lighthill - Whitham - Richards model:
\begin{equation}\label{LWR}
u_t + [u(1-u)]_x = 0.
\end{equation}
$\bullet$ (Look-A) Look ahead relaxation model:
 \begin{equation}\label{LAR}
u_t + [u(1-u) e^{-\bar{u}}]_x = 0.
\end{equation}
$\bullet$ (Look-AB) Look-ahead relaxation-behind intensification model:
 \begin{equation}\label{LAB}
u_t + [u(1-u) e^{-\bar{u} +\tilde{u}}]_x = 0.
\end{equation}
Here, the global terms associated with \emph{constant} interaction potentials(\eqref{traffic_vanilla_kernel} and \eqref{look_b_const}) are given by 
\begin{equation}\label{GT_CONST}
\bar{u}(t,x)=\frac{1}{\gamma_a}\int^{x+\gamma_a} _{x} u(t,y) \, dy, \ \ \ \ \ \tilde{u}(t,x):=\frac{1}{\gamma_b}\int^{x} _{x-\gamma_b} u(t,y) \, dy.
\end{equation}
Also, the linear interaction potentials(\eqref{linear_kernel} and \eqref{look_b_linear}) lead to
\begin{equation}\label{GT_LIN}
\bar{u}(t,x)=\frac{2}{\gamma_a}\int^{x+\gamma_a} _{x}  \bigg{(} 1+ \frac{x-y}{\gamma_a}  \bigg{)} u(t,y) \, dy, \ \ \ \ \ \tilde{u}(t,x)=\frac{2}{\gamma_b}\int^{x} _{x-\gamma_b}  \bigg{(} 1- \frac{x-y}{\gamma_b}  \bigg{)} u(t,y) \, dy.
\end{equation}

In all examples, we should point out that our look ahead-behind model is better in the sense of dispersive and smoothing effects.

\begin{figure}[H]\label{Example1}
\begin{subfigure}{.5\textwidth}
  \centering
  \includegraphics[width=1\linewidth]{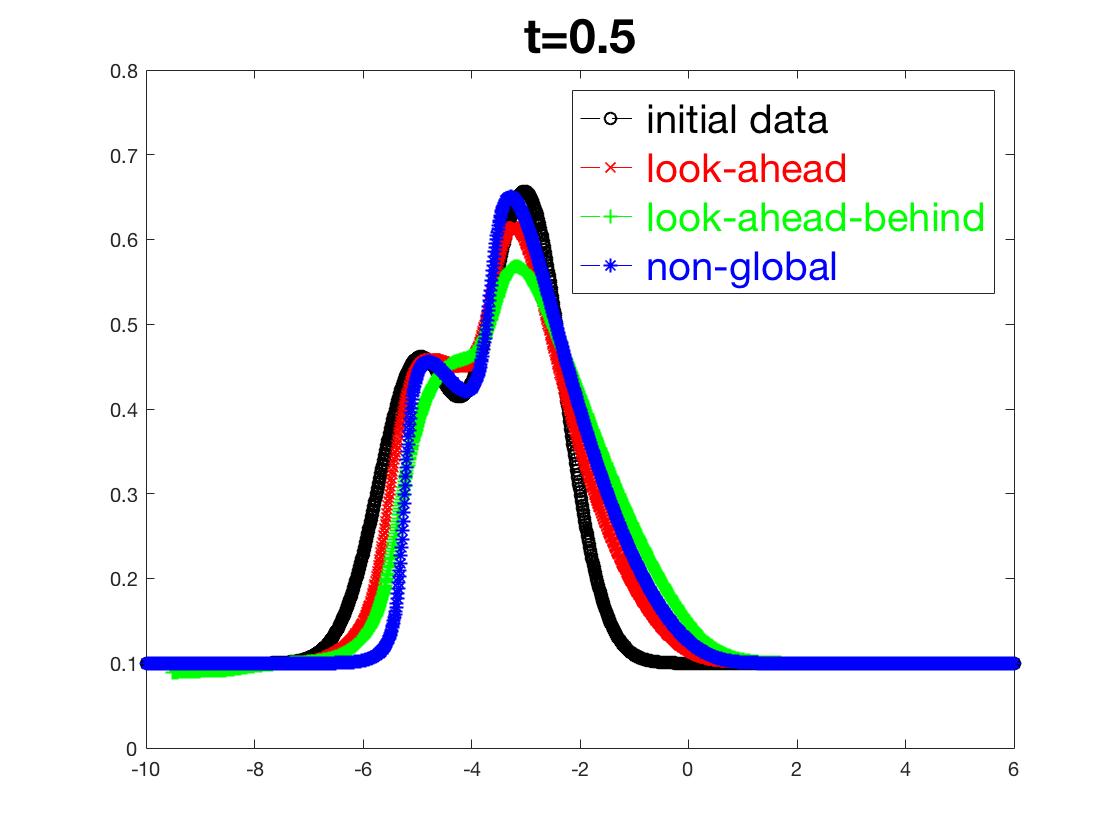}
  \label{fig:sfig1}
\end{subfigure}%
\begin{subfigure}{.5\textwidth}
  \centering
  \includegraphics[width=1\linewidth]{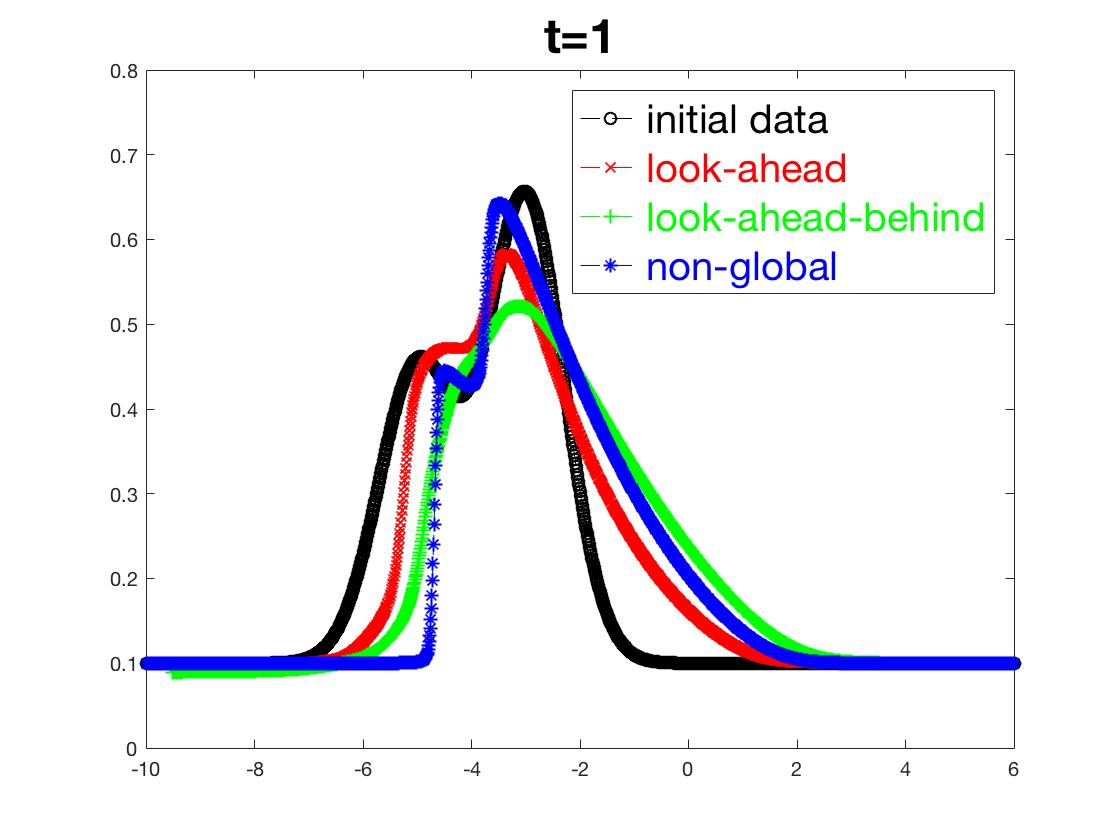}
  \label{fig:sfig2}
\end{subfigure}
\begin{subfigure}{.5\textwidth}
  \centering
  \includegraphics[width=1\linewidth]{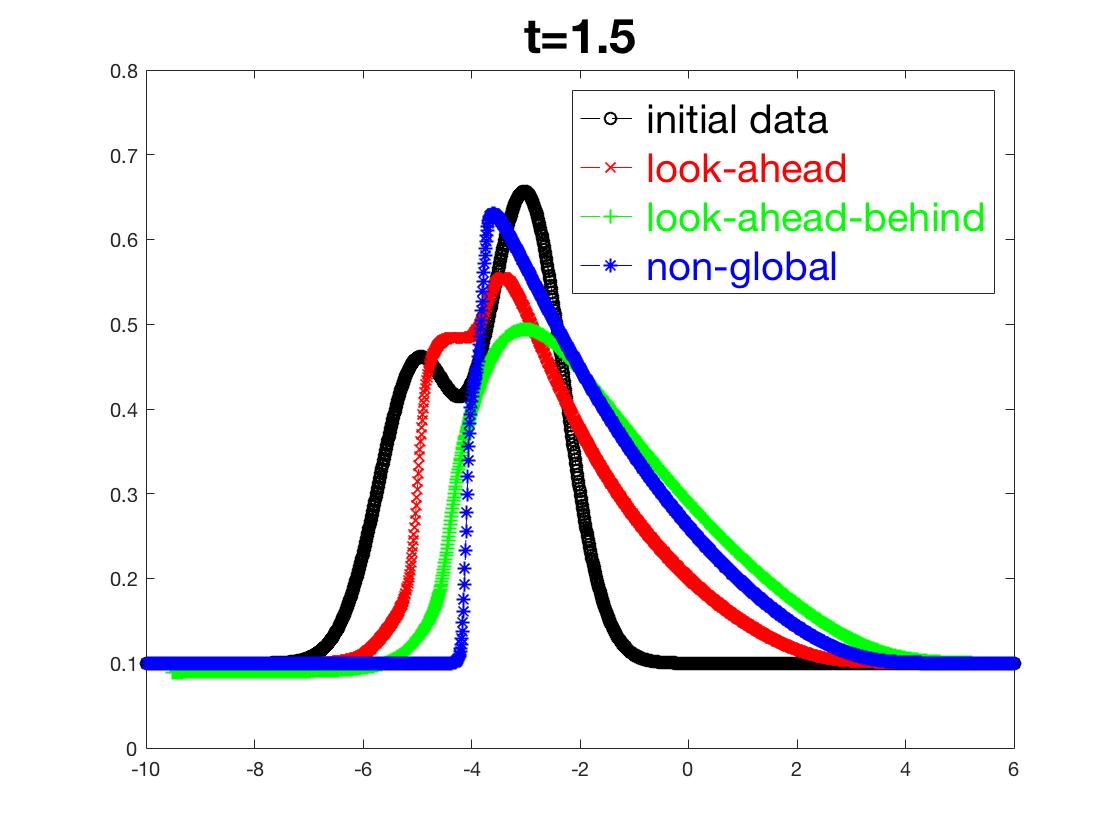}
  \label{fig:sfig1}
\end{subfigure}%
\begin{subfigure}{.5\textwidth}
  \centering
  \includegraphics[width=1\linewidth]{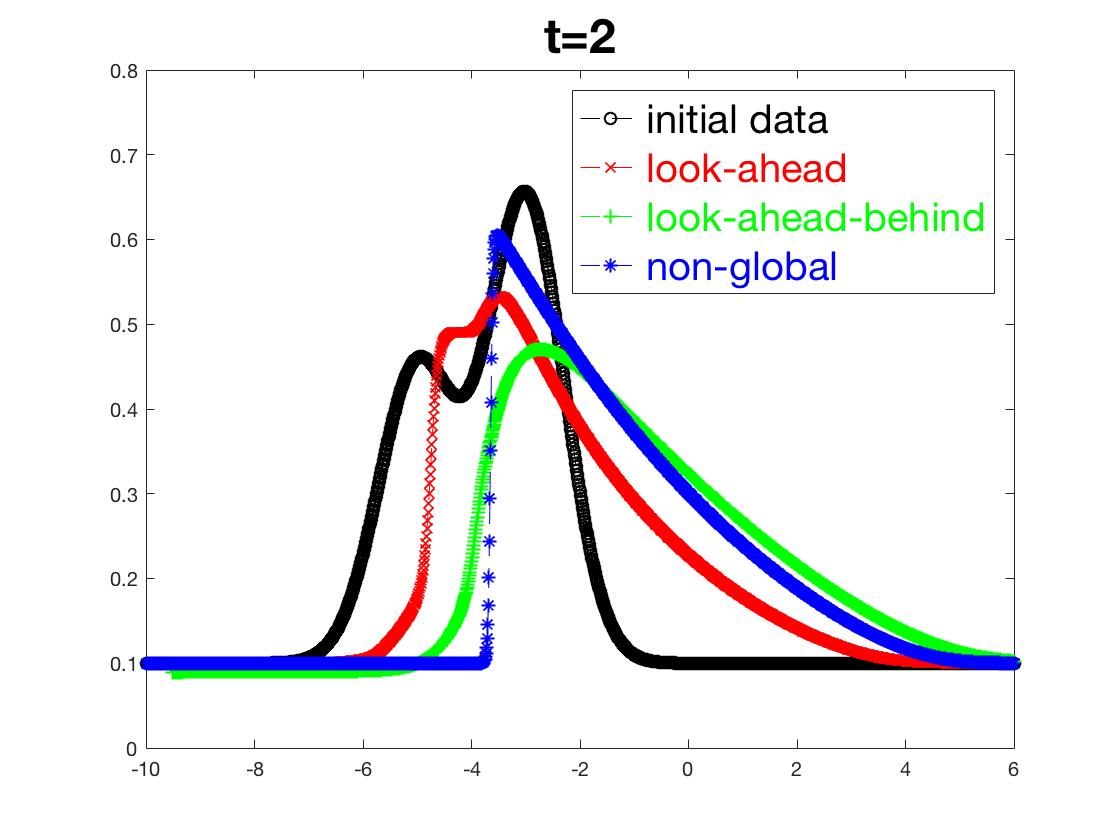}
  \label{fig:sfig2}
\end{subfigure}
\begin{subfigure}{.5\textwidth}
  \centering
  \includegraphics[width=1\linewidth]{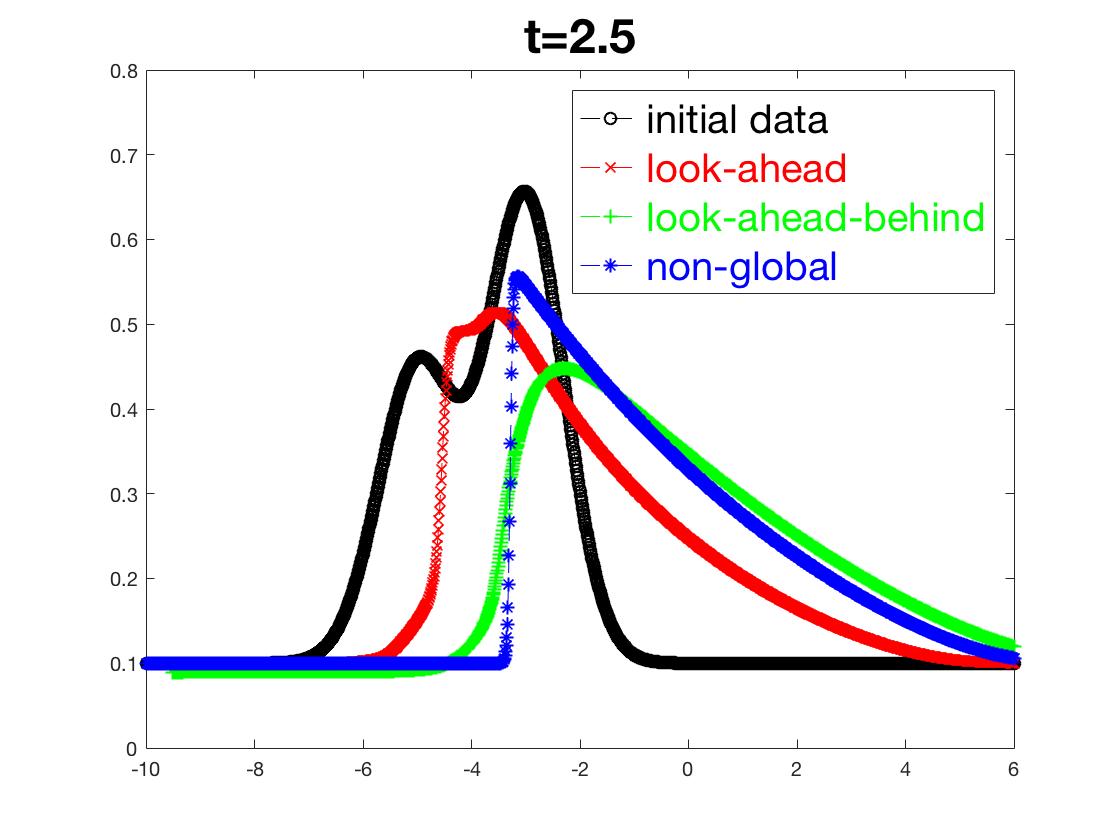}
  \label{fig:sfig1}
\end{subfigure}%
\begin{subfigure}{.5\textwidth}
  \centering
  \includegraphics[width=1\linewidth]{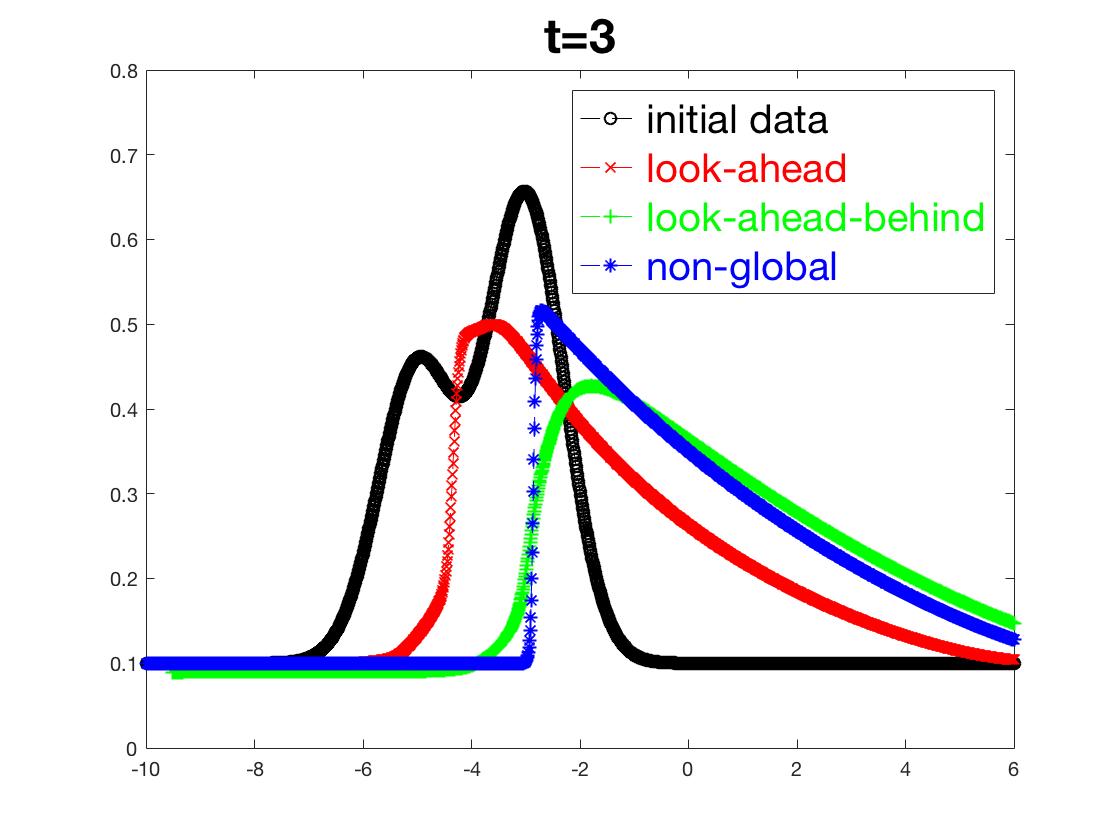}
  \label{fig:sfig2}
\end{subfigure}
\caption{Example \ref{eg1}: Time evolutions of the solutions equipped with LWR(Non-global flux), Look-Ahead flux and Ahead-Behind flux; $\gamma_a=1$ and $\gamma_b=0.5$.}
\label{fig_p1}
\end{figure}

\begin{example}\label{eg1}\emph{(two plateaus flow)}
We consider equations LWR, Look-A and Look-AB with constant interaction potentials(that is, \eqref{LWR}, \eqref{LAR}-\eqref{GT_CONST} and \eqref{LAB}-\eqref{GT_CONST}, respectively)
with $\gamma_a =1$ and $\gamma_b=0.5$ subject to the following initial data:
\begin{equation}
u(0,x)=0.1+0.35 e^{-(x+5)^2}+0.55 e^ {-(x+3)^2}.
\end{equation}
The situation corresponds to highly congested traffic($-4<x<2$) is preceded by free traffic and followed by less congested traffic. See Figure \ref{fig_p1}.
\end{example}

In this example, on can clearly see the effects of the global fluxes. Look-A model's waves(Red) lag behind that of the LWR waves(Blue). The waves in Look-AB model have faster pace than the ones in Look-A. This can be explained by the awareness of traffic behind. While the waves with global fluxes develop no shock, the LWR wave(Blue) develop a shock discontinuity at $t=2$.  It is interesting to observe that the valley between two initial plateaus is filled quickly in Look-AB model, while Look-A model has persistent valley until $t=2.5$. Overall, we should point out that Look-AB model has both enhanced dispersive and smoothing effects.



\begin{figure}[H]
\begin{subfigure}{.5\textwidth}
  \centering
  \includegraphics[width=1\linewidth]{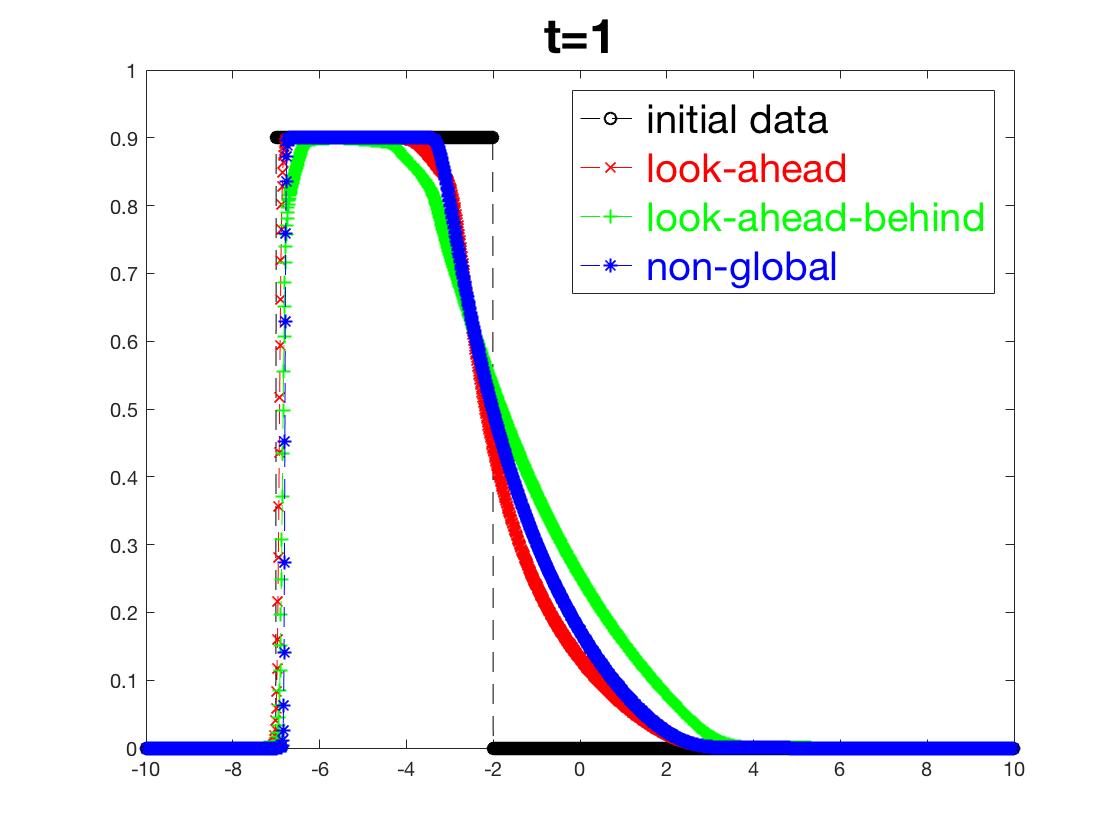}
  \label{fig:sfig1}
\end{subfigure}%
\begin{subfigure}{.5\textwidth}
  \centering
  \includegraphics[width=1\linewidth]{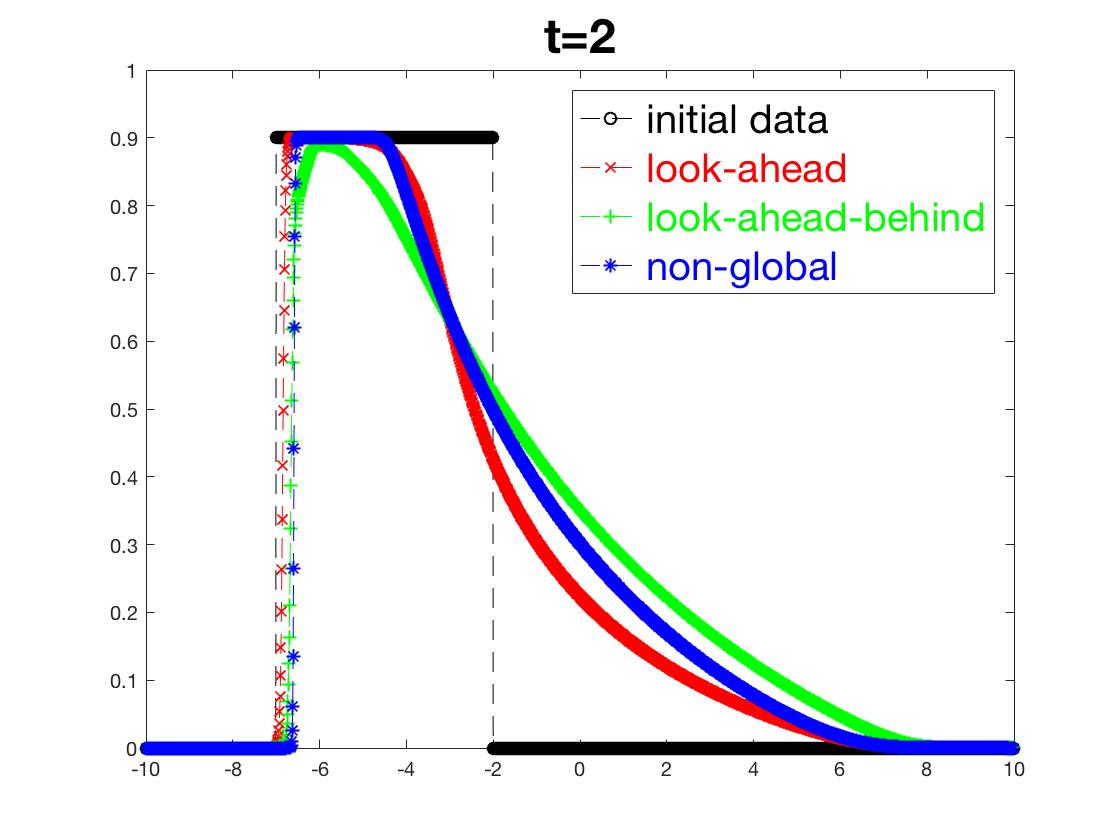}
  \label{fig:sfig2}
\end{subfigure}
\begin{subfigure}{.5\textwidth}
  \centering
  \includegraphics[width=1\linewidth]{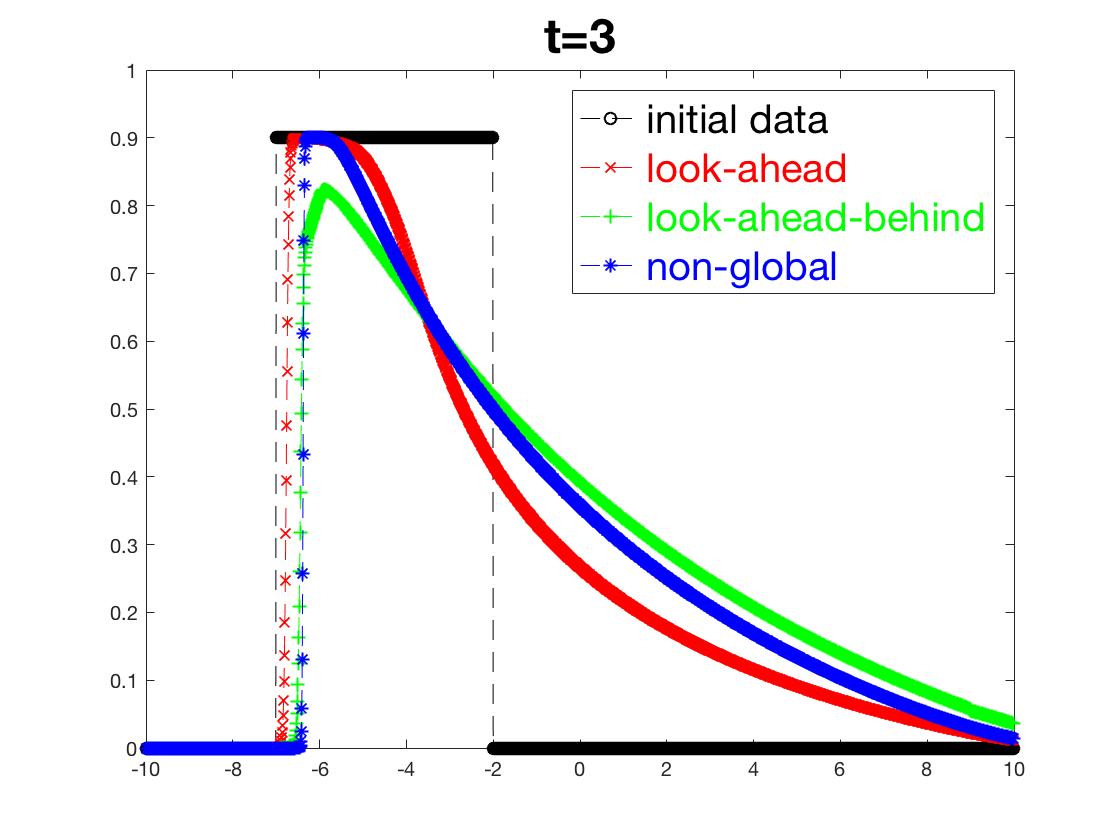}
  \label{fig:sfig1}
\end{subfigure}%
\begin{subfigure}{.5\textwidth}
  \centering
  \includegraphics[width=1\linewidth]{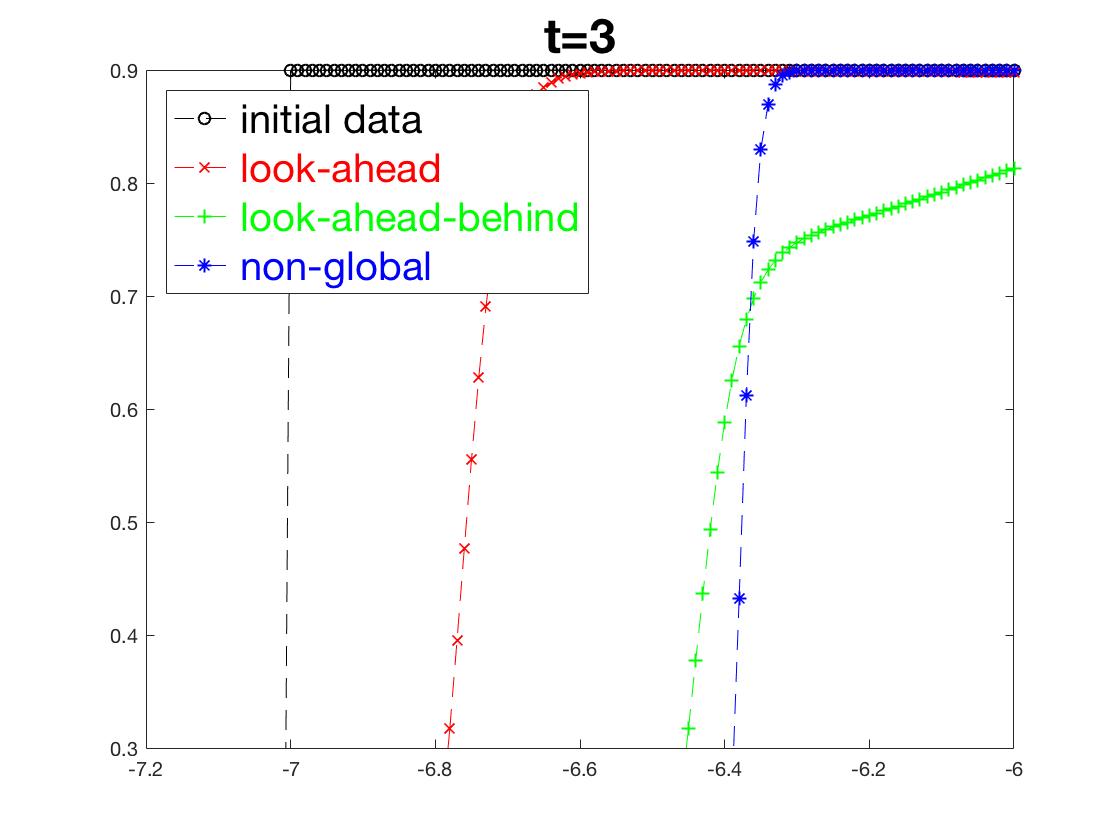}
  \label{fig:sfig2}
\end{subfigure}
\caption{Example \ref{eg2}: Time evolutions of the solutions equipped with LWR(Non-global flux), Look-Ahead flux and Ahead-Behind flux; $\gamma_a=1$ and $\gamma_b=0.5$, and zoom near the steepest gradients(lower right).}
\label{fig_p2}
\end{figure}

\begin{example}\label{eg2}\emph{(red light
 traffic)}
 We consider equations LWR, Look-A and Look-AB with constant interaction potentials(that is, \eqref{LWR}, \eqref{LAR}-\eqref{GT_CONST} and \eqref{LAB}-\eqref{GT_CONST}, respectively) with $\gamma_a =1$ and $\gamma_b=0.5$ subject to the following initial data:
\begin{equation}
u(0,x)=
\left\{
  \begin{array}{ll}
0.9, & x \in [-7,-2], \\
 0, &   \mathrm{otherwise}.
  \end{array}
\right.
\end{equation}
This initial condition corresponds to a situation in which the red traffic light is located at $x=-2$ and it is turned to green at the initial time instance.  See Figure \ref{fig_p2}.
\end{example}

In this example, it is observed that waves in the Look-AB model have the fastest pace and the most prominent dispersive effect. Furthermore, if we zoom in near the steepest gradients(lower right), one can clearly see that the Look-AB wave has less steep slope. It looks like the behind intensification dynamics in the Look-AB model not only increase the flux, but also contribute to the dispersive effect through global interactions.

\begin{figure}[H]
\begin{subfigure}{.33\textwidth}
  \centering
  \includegraphics[width=1\linewidth]{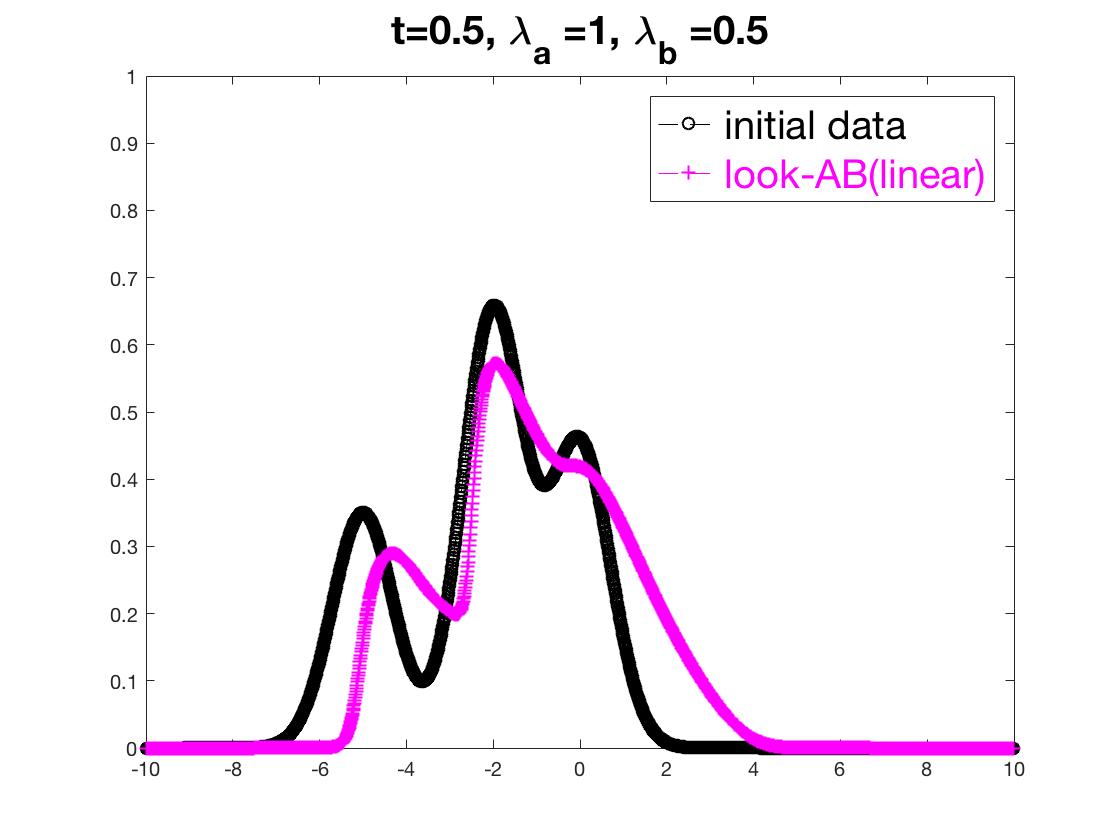}
  \label{fig:sfig1}
\end{subfigure}%
\begin{subfigure}{.33\textwidth}
  \centering
  \includegraphics[width=1\linewidth]{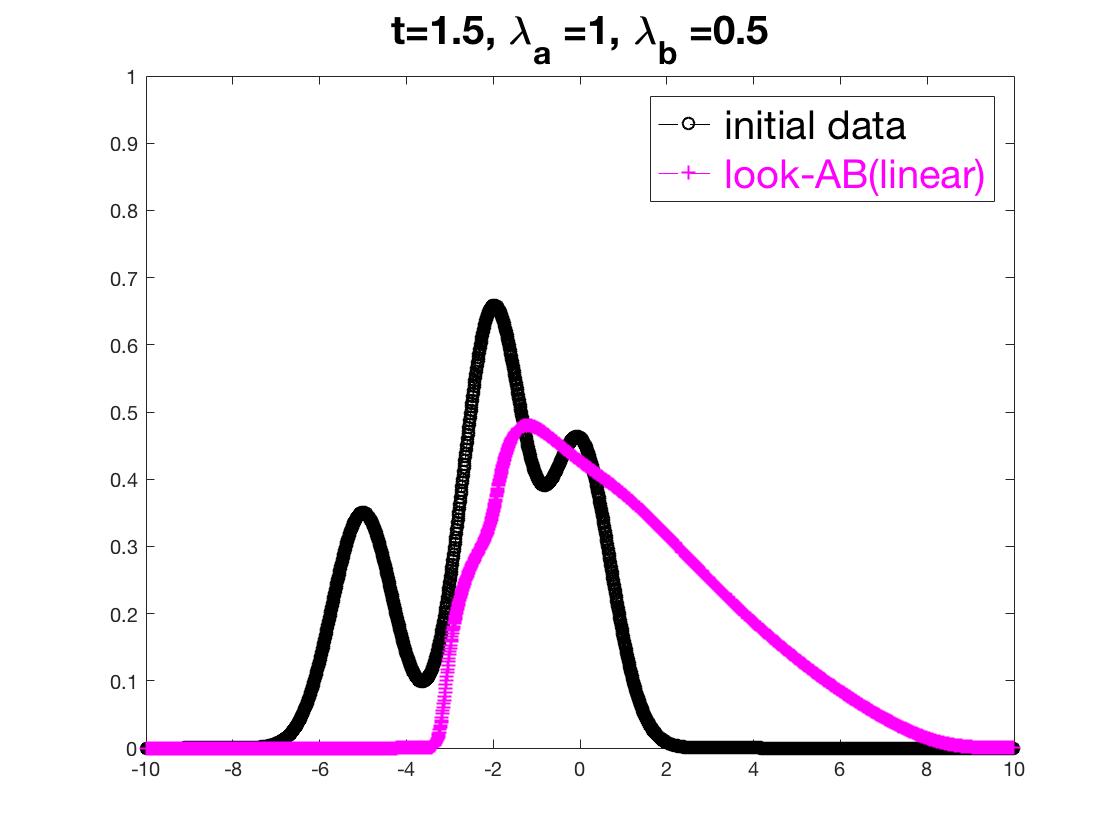}
  \label{fig:sfig2}
\end{subfigure}
\begin{subfigure}{.33\textwidth}
  \centering
  \includegraphics[width=1\linewidth]{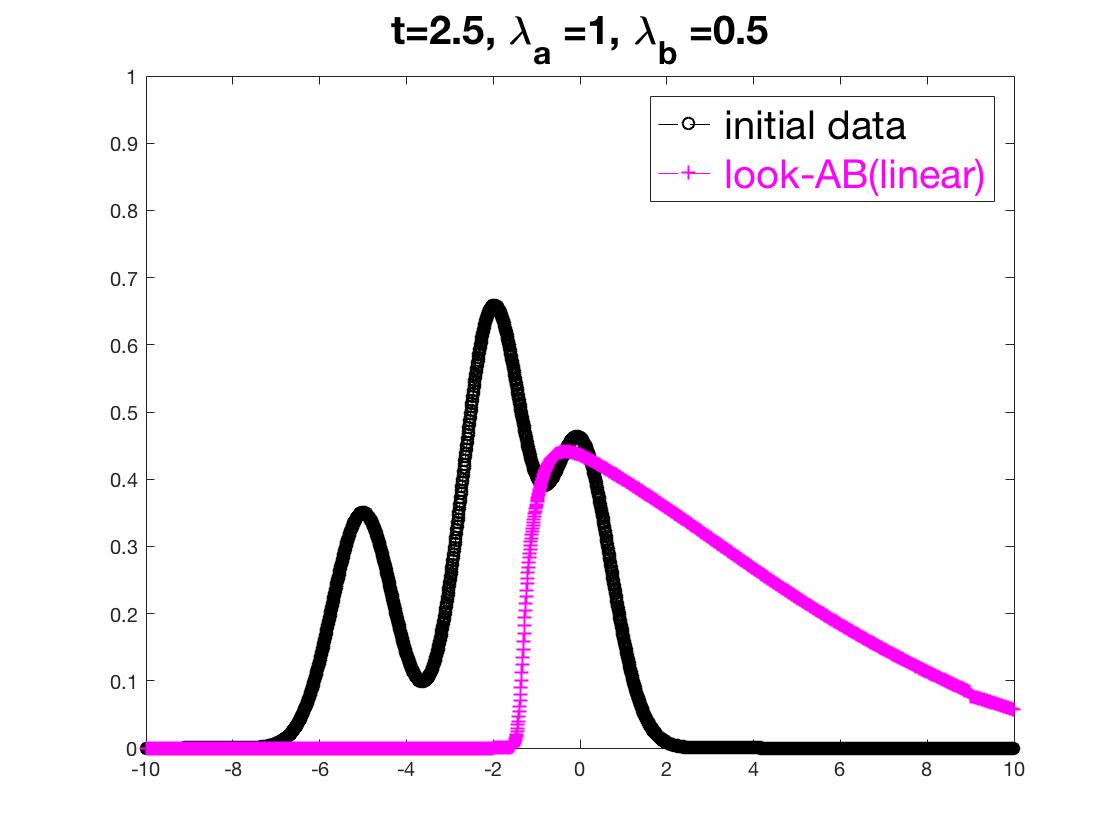}
  \label{fig:sfig1}
\end{subfigure}\\

\begin{subfigure}{.33\textwidth}
  \centering
  \includegraphics[width=1\linewidth]{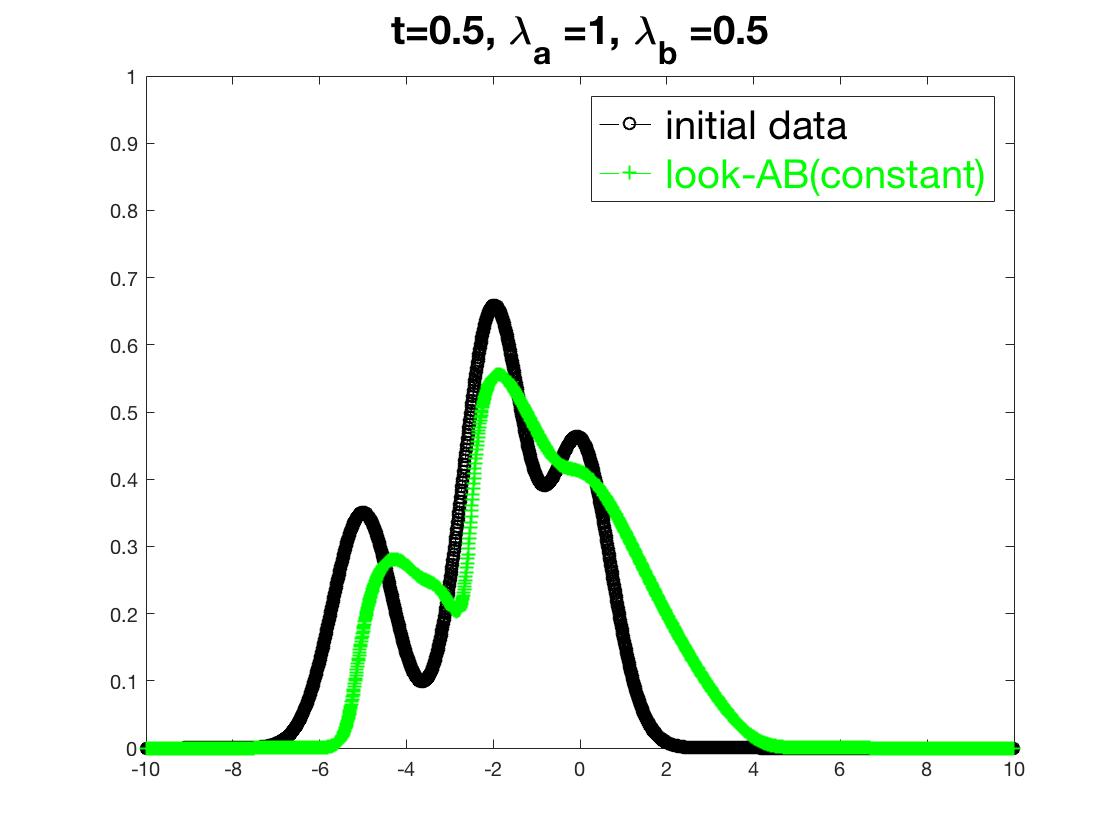}
  \label{fig:sfig1}
\end{subfigure}%
\begin{subfigure}{.33\textwidth}
  \centering
  \includegraphics[width=1\linewidth]{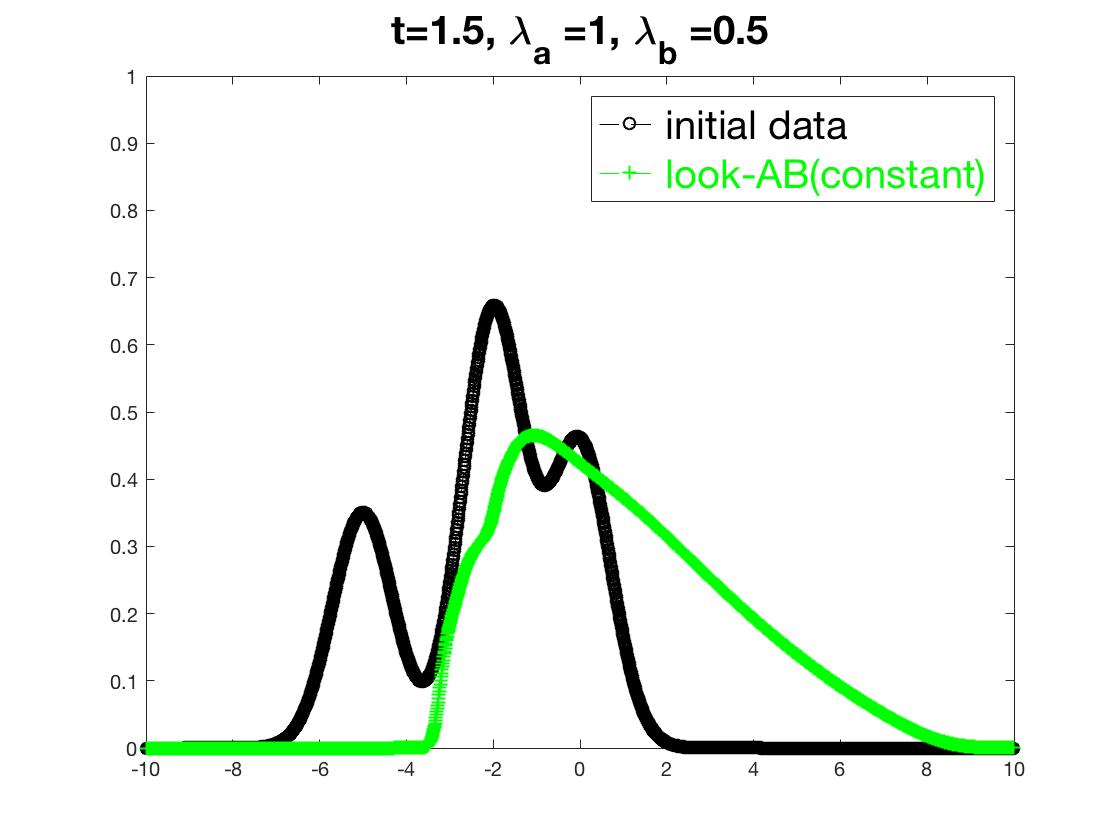}
  \label{fig:sfig2}
\end{subfigure}
\begin{subfigure}{.33\textwidth}
  \centering
  \includegraphics[width=1\linewidth]{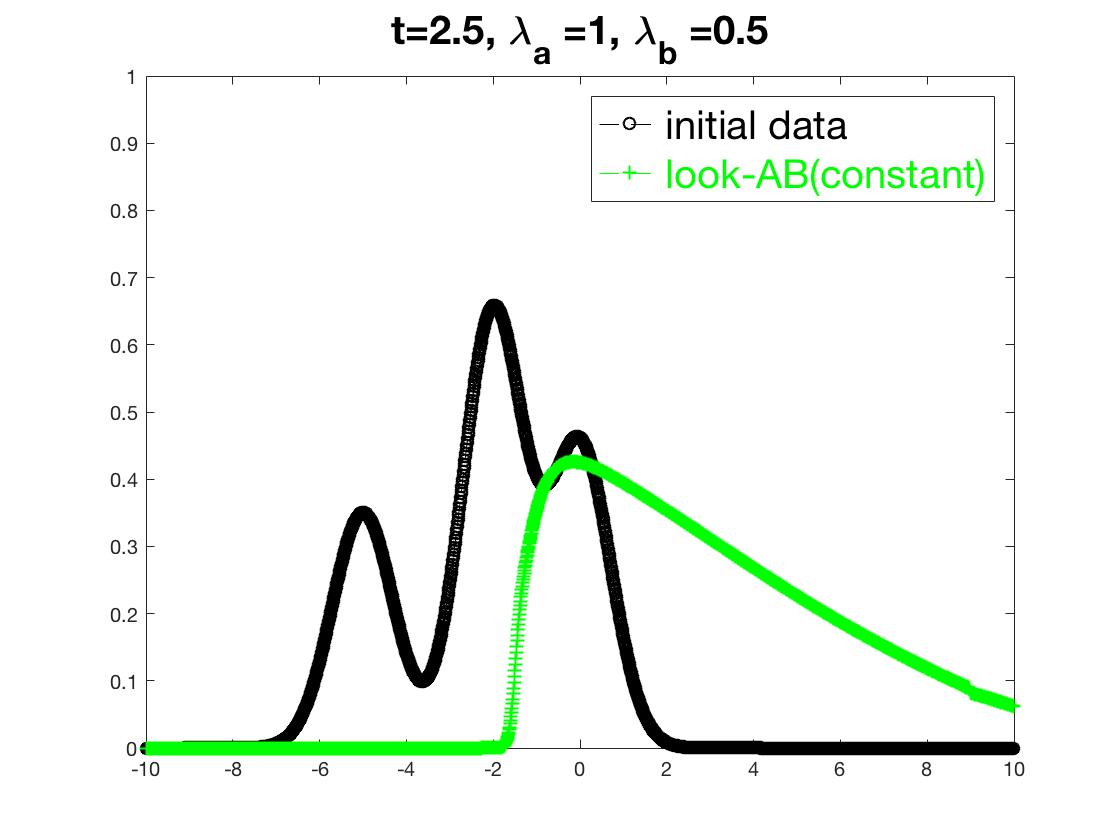}
  \label{fig:sfig1}
\end{subfigure}\\

\begin{subfigure}{.33\textwidth}
  \centering
  \includegraphics[width=1\linewidth]{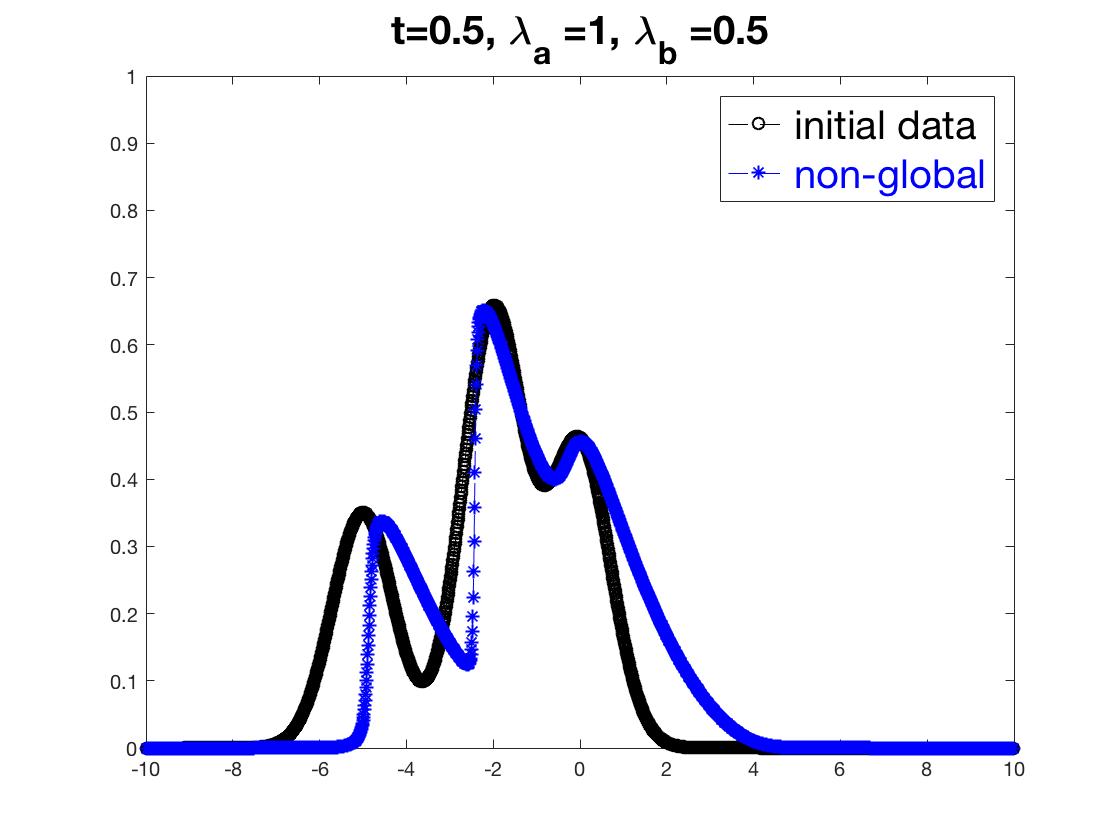}
  \label{fig:sfig1}
\end{subfigure}%
\begin{subfigure}{.33\textwidth}
  \centering
  \includegraphics[width=1\linewidth]{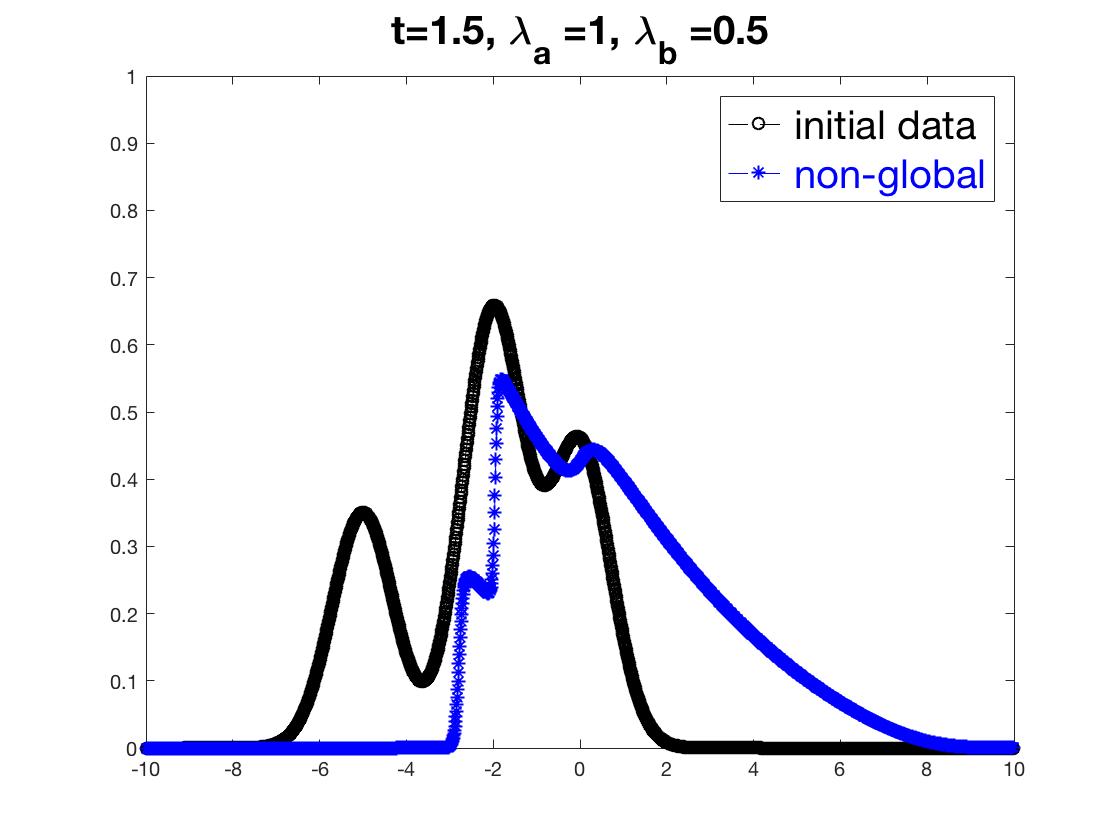}
  \label{fig:sfig2}
\end{subfigure}
\begin{subfigure}{.33\textwidth}
  \centering
  \includegraphics[width=1\linewidth]{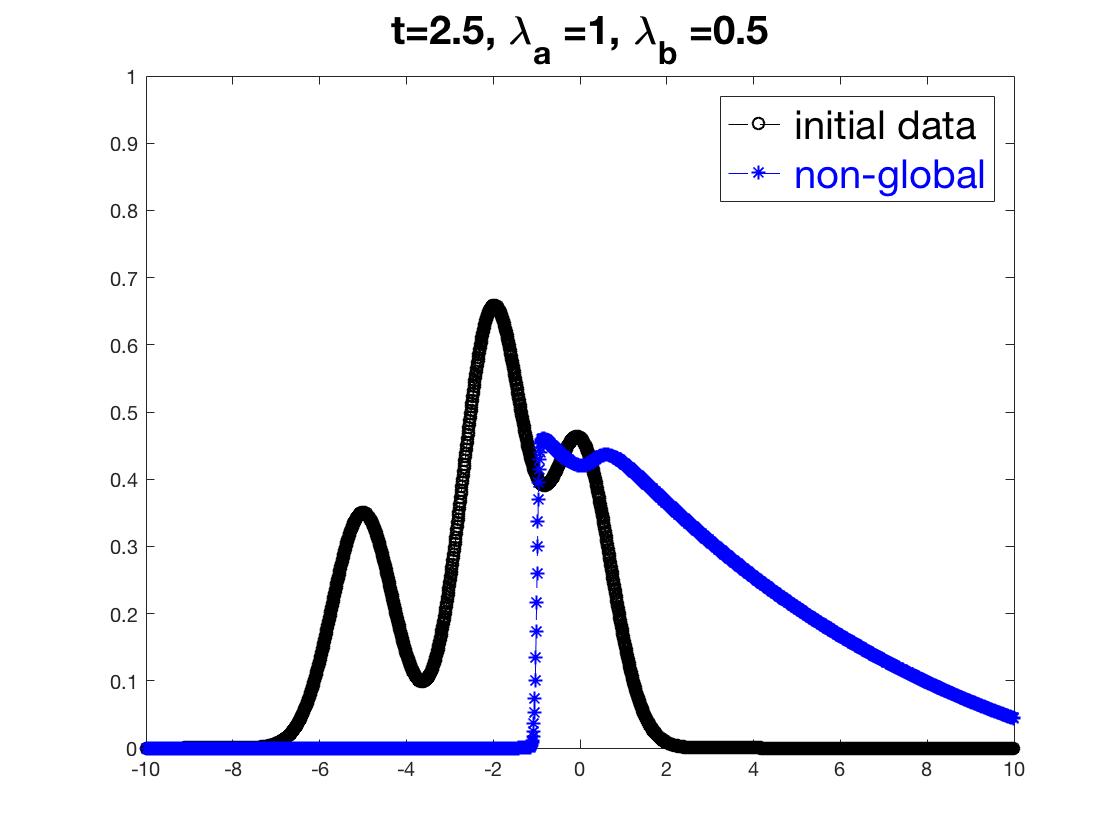}
  \label{fig:sfig1}
\end{subfigure}\\

\caption{ Example \ref{eg3}: Time evolutions(left to right) of Look-AB with linear interaction kernel(upper) and constant interaction kernel(middle), and LWR(lower).}
\label{fig_p51}
\end{figure}

\begin{example}\label{eg3}\emph{(three plateaus; constant interaction potentials and linear interaction potentials)}.
We consider Look-AB with constant interaction potentials and linear interaction potentials(\eqref{LAB}-\eqref{GT_LIN} and \eqref{LAB}-\eqref{GT_CONST}), subject to the following initial data
\begin{equation}
u(0,x)=0.35 e^{-(x+5)^2}+0.65 e^ {-(x+2)^2} + 0.45 e^{-x^2},
\end{equation}
with $ (\gamma_a,  \gamma_b) = (1, 0.5)$.  See Figures \ref{fig_p51}.
\end{example}
In this example, we observe that the LWR model's solution(blue) develops the shock. In contrast to this, the Look-AB solutions with linear potentials(magenta) and constant potentials(green) have no shock. Furthermore, the valleys between three plateaus are filled relatively quickly in Look-AB models. 


\begin{figure}[H]
\begin{subfigure}{.33\textwidth}
  \centering
  \includegraphics[width=1\linewidth]{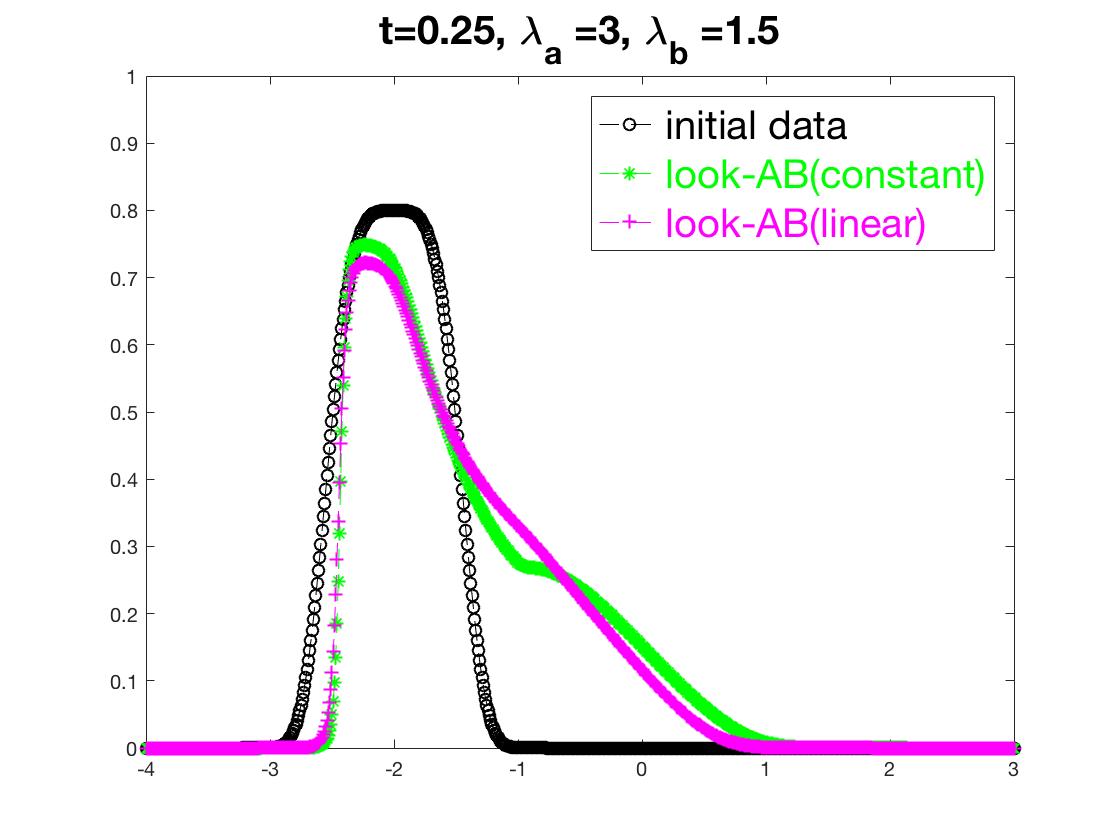}
  \label{fig:sfig1}
\end{subfigure}%
\begin{subfigure}{.33\textwidth}
  \centering
  \includegraphics[width=1\linewidth]{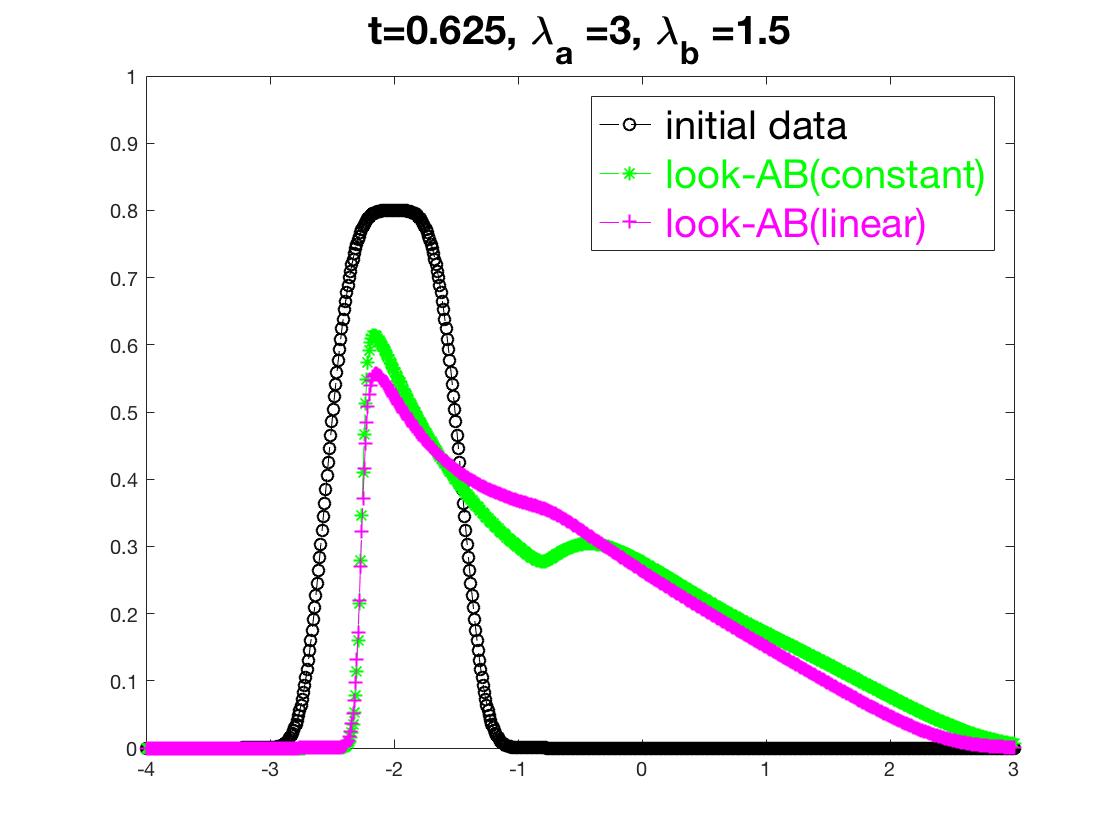}
  \label{fig:sfig2}
\end{subfigure}
\begin{subfigure}{.33\textwidth}
  \centering
  \includegraphics[width=1\linewidth]{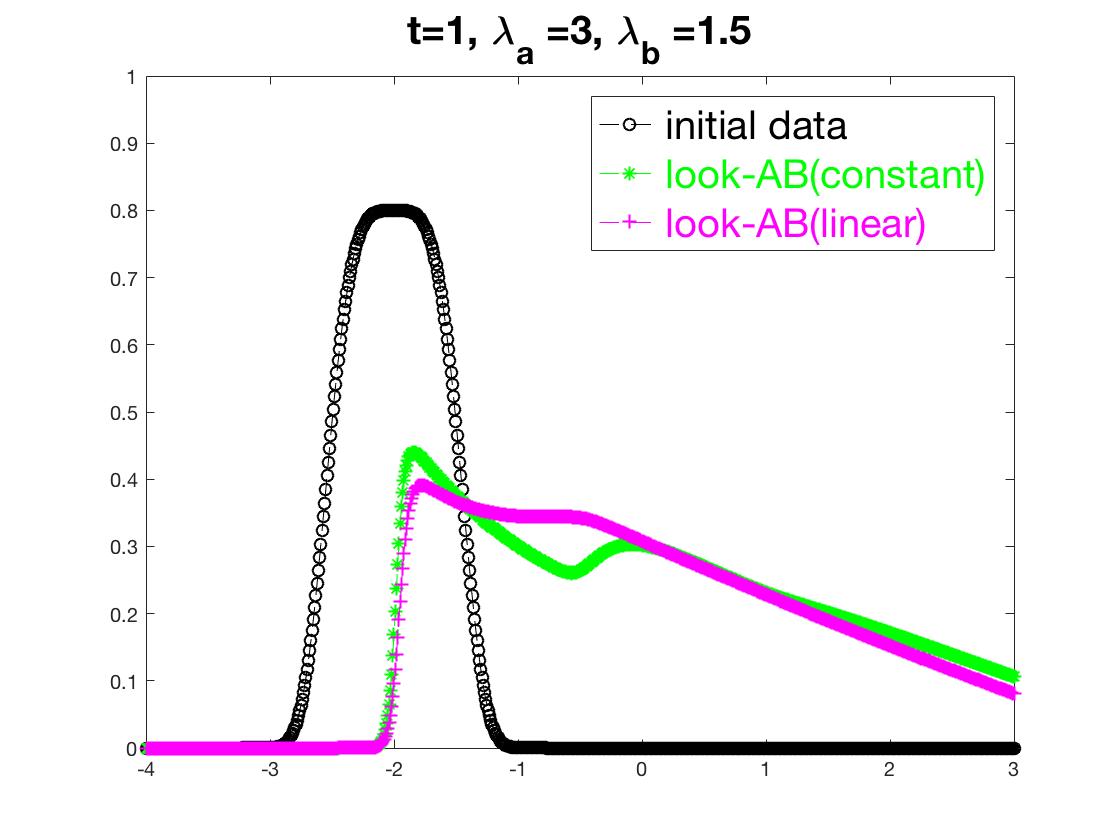}
  \label{fig:sfig1}
\end{subfigure}\\

\begin{subfigure}{.33\textwidth}
  \centering
  \includegraphics[width=1\linewidth]{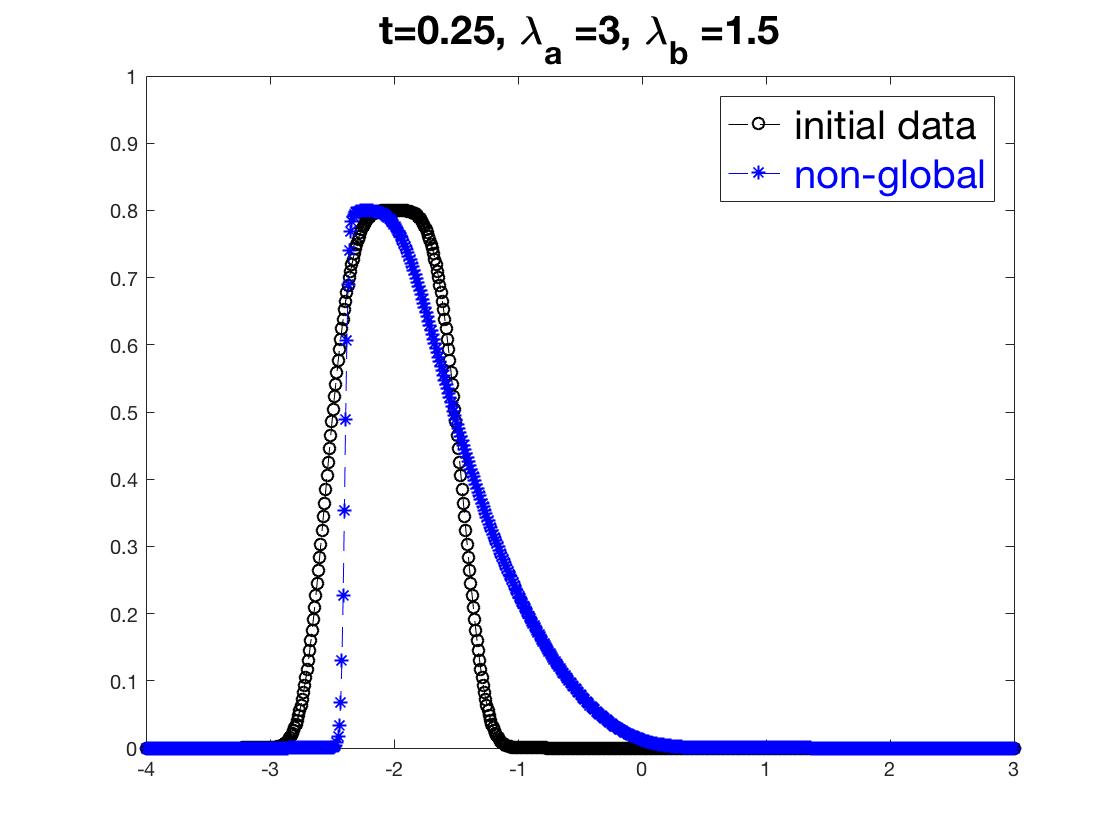}
  \label{fig:sfig1}
\end{subfigure}%
\begin{subfigure}{.33\textwidth}
  \centering
  \includegraphics[width=1\linewidth]{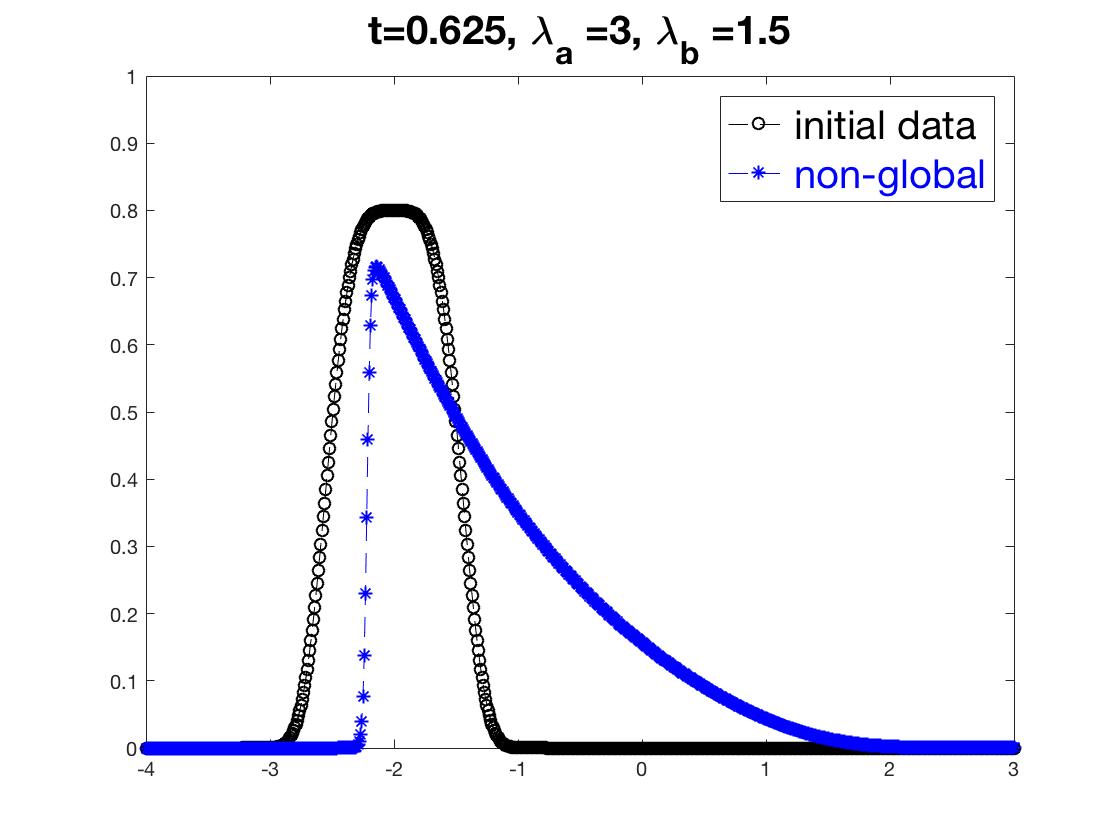}
  \label{fig:sfig2}
\end{subfigure}
\begin{subfigure}{.33\textwidth}
  \centering
  \includegraphics[width=1\linewidth]{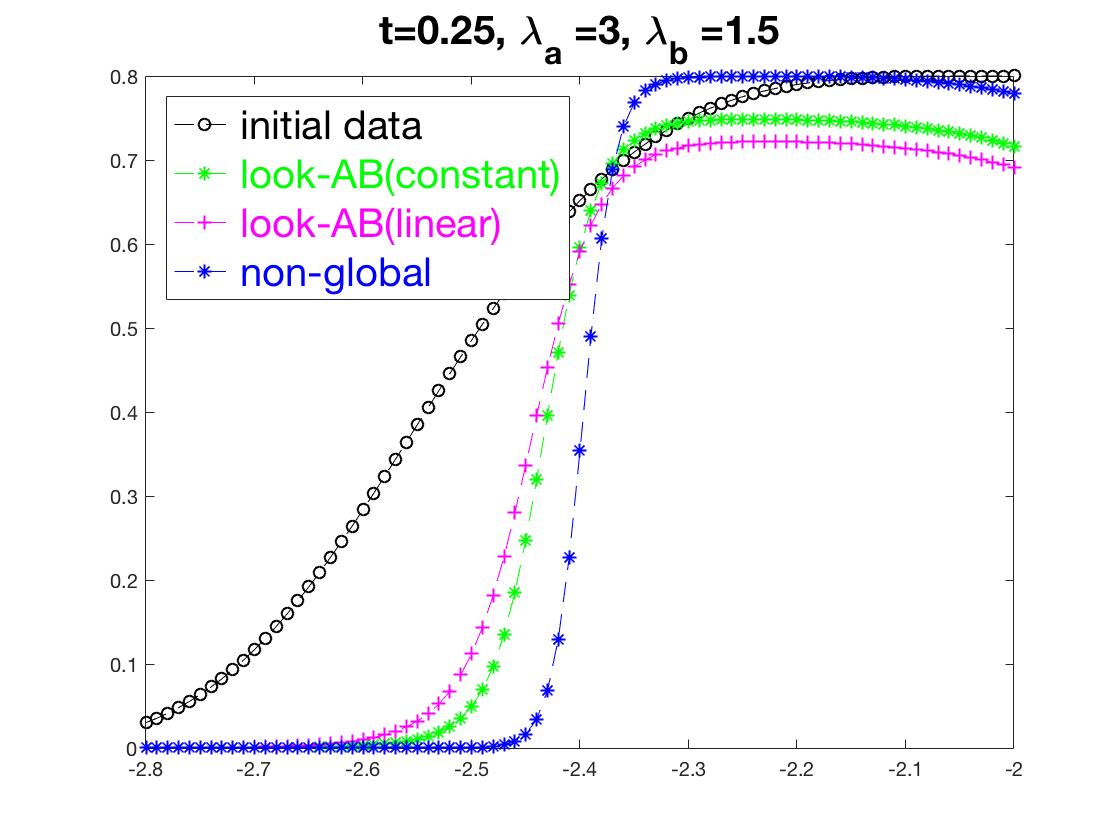}
  \label{fig:sfig1}
\end{subfigure}
\caption{ Example \ref{eg4} Time evolutions(left to right) of the solutions equipped with; linear interaction kernel(magenta), constant interaction kernel(green) and LWR(blue), zoom in near the steepest slopes(lower right).}
\label{fig_p52}
\end{figure}

\begin{example}\label{eg4}\emph{(steep plateau; constant interaction potentials and linear interaction potentials)}
We consider Look-AB with constant interaction potentials and linear interaction potentials(\eqref{LAB}-\eqref{GT_LIN} and \eqref{LAB}-\eqref{GT_CONST}), subject to the following initial data
\begin{equation}
u(0,x)=0.80 e^{-8(x+2)^4},
\end{equation}
with $ (\gamma_a,  \gamma_b) = (3, 1.5)$.  See Figures \ref{fig_p52}.
\end{example}

With this very steep initial wave,  one can again observe that the waves in Look-AB model have enhance dispersive effects compare to the LWR wave(blue). One can clearly see the shock formation in the LWR wave. Furthermore, it is observed that the Look-AB model with linear interaction potentials(magenta) has less steep wave than the model with constant potentials(green); see the lower right figure in Figure \ref{fig_p52}. This observation is consistent with remark iii) in the introduction.

\bibliographystyle{abbrv}

\begin{thebibliography}{10}

\bibitem{BBKT}
F.~Betancourt, R.~Burger, K.~H.~Karlsen and E.~M.~ Tory.
\newblock  On nonlocal conservation laws modelling sedimentation.
\newblock {\em Nonlinearity.},  24: 855--885,  2011.

\bibitem{MDS08}
M. Burger and Y. Dolak and C. Schmeiser.
\newblock Asymptotic analysis of an advection-dominated chemotaxis model in multiple spatial dimensions.
\newblock {\em Commun.  Math. Sci. }, 6(1):1--28, 2008.


\bibitem{CE98}
A. Constantin and J. Escher.
\newblock Wave breaking for nonlinear nonlocal shallow water equations.
\newblock  {\em Acta Math.}, 181: 229--243, 1998.

 \bibitem{DS05}
Y. Dolak and C. Schmeiser.
\newblock The Keller-Segel model with logistic sensitivity function and small diffusivity.
\newblock  {\em SIAM J. Appl. Math.}, 66: 286--308, 2005.

\bibitem{ELT01}
S. Engelberg,  H. Liu and E. Tadmor
\newblock  Critical Thresholds in Euler-Poisson equations.
\newblock {\em Indiana Univ. Math. J.},  50: 109--157,  2001.



\bibitem{Holm}
D. D. ~Holm and A. N. W. ~Hone.
\newblock  A class of equations with peakon and pulson solutions (with a appendix by Braden H and Byatt-Smith).
\newblock {\em J. Nonlinear Math. Phys},  12 Suppl. 1: 380--394,  2005.

\bibitem{Hunter}
J. K. ~Hunter.
\newblock  Numerical solutions of some nonlinear dispersive wave equations.
\newblock {\em Lect. Appl. Math},  301--316,  1990.


\bibitem{KP09}
A. Kurganov and A. Polizzi.
\newblock Non-oscillatory central schemes for traffic flow models with Arrhenius look-ahead dynamics.
\newblock  {\em Netw. Heterog. Media.},   4: 431--451, 2009.

\bibitem{Kynch}
G.~Kynch.
\newblock  A theory of sedimentation.
\newblock {\em Trans. Fraday Soc.},  48: 66--76,  1952.




\bibitem{LL15}
Y. ~Lee and  H. Liu.
\newblock Thresholds for shock formation in traffic flow models with Arrhenius look-ahead dynamics.
\newblock  {\em  DCDS-A},  35(1): 323--3339, 2015.


\bibitem{LL11}
D. Li, and T. Li.
\newblock Shock formation in a traffic flow model with Arrhenius look-ahead dynamics.
\newblock  {\em Netw. Heterog. Media.}, 6: 681--694, 2011


\bibitem{LL09}
T. Li, and H. Liu.
\newblock Critical thresholds in hyperbolic relaxation systems.
\newblock  {\em J. Differential Equations}, 247: 33--48, 2009




\bibitem{LW55}
M. J. Lighthill and G. B. Whitham. 
\newblock On kinematic waves: II. A theory of traffic flow on long crowded roads.
\newblock  {\em  Proc. Roy. Soc., London, Ser. A}, 229: 317--345, 1955

\bibitem{Liu0}
H.~Liu.
\newblock  Wave breaking in a class of nonlocal dispersive wave equations.
\newblock {\em Journal of Nonlinear Math Phys.},  13(3): 441--466,  2006.

\bibitem{LT02}
H. ~Liu, E. Tadmor.
\newblock Spectral dynamics of the velocity gradient field in restricted fluid flows.
 \newblock  {\em Comm. Math. Phys.},   228: 435--466,  2002.



\bibitem{Vak1}
E. J. ~Parkes and V. O. ~Vakhneko.
\newblock  The calculation of multi-soliton solutions of the Vakhnenko equation by the inverse scattering method.
\newblock {\em Chaos Solitons Fractals.},  13: 1819--1826,  2002.



\bibitem{R56}
P. I. Richards. 
\newblock Shock waves on the highway.
\newblock  {\em Oper. Res}, 4: 42--51, 1956

 \bibitem{JR90}
J. Rubinstein.
\newblock Evolution equations for stratified dilute suspensions.
\newblock  {\em Phys. Fluids A.}, 2(1): 3--6, 1990.

 \bibitem{RK89}
J. Rubinstein and J. B. Keller.
\newblock Sedimentation of a dilute suspension.
\newblock  {\em Phys. Fluids A.}, 1: 637--643, 1989.


\bibitem{SK06}
A. Sopasakis and M. Katsoulakis.
\newblock Stochastic modeling and simulation of traffic flow: Asymmetric single exclusion process with Arrhenius look-ahead dynamics.
\newblock {\em SIAM J. Appl. Math.},  66(3): 921--944, 2006.

\bibitem{TT14}
E. Tadmor and C. Tan
\newblock Critical thresholds in flocking hydrodynamics with non-local alignment.
\newblock  {\em  Philos Trans A Math Phys Eng Sc.},  372(2028), 2014.


\bibitem{GW74}
G. B. Whitham.
\newblock Linear and nonlinear waves.
\newblock  {\em John Wiley and Sons}, 1974

 \bibitem{KZ99}
K. Zumbrun
\newblock On a nonlocal dispersive equation modeling particle suspensions.
\newblock  {\em Quart. Appl. Math.}, 57: 573--600, 1999.



\end{thebibliography}

\end{document}